\documentclass{conm-p-l}

\usepackage{amsmath}
\usepackage{amscd}
\usepackage{amssymb}

\newtheorem{theorem}{Theorem}[section]

\theoremstyle{definition}
\newtheorem{corollary}[theorem]{Corollary}
\newtheorem{definition}[theorem]{Definition}

\newtheorem{proposition}[theorem]{Proposition}
\newtheorem{remark}[theorem]{Remark}

\theoremstyle{remark}

\newcommand{\ad}{\text{ad}}

\newcommand{\ch}{\mbox{ch} }
\newcommand{\C}{ \mathbb C }

\newcommand{\End}{{\rm End}}

\newcommand{\Fock}{{\mathcal F}_X}
\newcommand{\Fockh}{{\mathcal H}_X}
\newcommand{\g}{{\gamma}}
\newcommand{\G}{{\Gamma}}
\newcommand{\Gn}{{\G_n}}
\newcommand{\Gh}{{\mathfrak G}_{\hbar}}
\newcommand{\Hxn}{H^*(\Xn)}

\newcommand{\hf}{\frac12}
\newcommand{\Ind}{\mbox{Ind}}
\newcommand{\jj}{\mathfrak J}
\newcommand{\KX}{K}
\newcommand{\la}{\lambda}
\newcommand{\lambsq}{s(\lambda)}
\newcommand{\Oh}{{\mathfrak O}_{\hbar}}
\newcommand{\orbsym}{H^*_{\text{orb}}(X^n/S_n)}

\newcommand{\Res}{ {\rm Res} }
\newcommand{\pa}{\partial}
\newcommand{\ta}{\tau_*\alpha}
\newcommand{\tea}{\tau_*(e\alpha)}

\newcommand{\vac}{|0\rangle}
\newcommand{\W}{{\mathcal W}}
\newcommand{\Wax}{{\mathcal W}_X}
\newcommand{\Xn}{ X^{[n]}}
\newcommand{\Z}{ \mathbb Z }

\numberwithin{equation}{section}

\begin{document}

\title[Hilbert schemes and symmetric products: a dictionary]
    {Hilbert schemes and symmetric products: \\a dictionary}

\author{Zhenbo Qin}
\address{Department of Mathematics, University of Missouri, Columbia, MO
65211, USA} \email{zq@math.missouri.edu}
\thanks{Both authors are partially supported by NSF grants}

\author{Weiqiang Wang}
\address{Department of Mathematics, University of Virginia,
Charlottesville, VA 22904} \email{ww9c@virginia.edu}


\begin{abstract}
Given a closed complex manifold $X$ of even dimension, we develop
a systematic (vertex) algebraic approach to study the rational
orbifold cohomology rings $\orbsym$ of the symmetric products. We
present constructions and establish results on the rings $\orbsym$
including two sets of ring generators, universality and stability,
as well as connections with vertex operators and $\mathcal W$
algebras. These are independent of but parallel to the main
results on the cohomology rings of the Hilbert schemes of points
on surfaces as developed in our earlier works joint with W.-P.~Li.
We introduce a deformation of the orbifold cup product and explain
how it is reflected in terms of modification of vertex operators
in the symmetric product case. As a corollary, we obtain a new
proof of the isomorphism between the rational cohomology ring of
Hilbert schemes $\Xn$ and the ring $\orbsym$ (after some
modification of signs), when $X$ is a projective surface with a
numerically trivial canonical class; we show that no sign
modification is needed if both cohomology rings use
$\C$-coefficients.
\end{abstract}
\maketitle

\date{}
\tableofcontents
\section{Introduction}

One of the recent surprises is the deep connections among geometry
of Hilbert schemes $\Xn$ of points on a projective surface $X$,
symmetric groups (or more generally wreath products), infinite
dimensional Lie algebras, and vertex algebras. The construction in
\cite{Na1, Gro} of Heisenberg algebra which acts on the direct sum
$\Fockh=\oplus_{n=0}^\infty \Hxn$ of rational cohomology groups of
Hilbert schemes made it possible to study the geometry of the
Hilbert schemes from a new algebraic viewpoint. Lehn's work
\cite{Lehn} initiated the deep interaction between Heisenberg
algebra and the ring structure on $H^*(\Xn)$. The potential of
this approach is made manifest in \cite{LQW1, LQW2} where two
different sets of ring generators on $H^*(\Xn)$ were found for an
arbitrary $X$. In \cite{LS2}, a construction of the ring $\Hxn$ is
made in terms of symmetric group, for $X$ with numerically trivial
canonical class. In \cite{LQW3}, the universality and stability of
Hilbert schemes were established, which concern about the
relations among the cohomology rings $\Hxn$ when $X$ or $n$
varies. In \cite{LQW4}, a $\W$ algebra was constructed
geometrically acting on $\Fockh$. This is an analogue of the
$\W_{1+\infty}$ algebra (cf. e.g. \cite{FKRW, Kac}) and it
contains the Heisenberg algebra and Virasoro algebra as
subalgebras. The construction of the $\W$ algebra is based on an
explicit vertex operator formula for the so-called Chern character
operator which plays an important role in the development.

It is well known that the Hilbert scheme $\Xn$ is a crepant
resolution of singularities of the symmetric product $X^n/S_n$. As
inspired from orbifold string theory \cite{DHVW, VW}, the geometry
of Hilbert schemes should be ``equivalent" to the
$S_n$-equivariant geometry of $X^n$. As stated in the
footnote~3 of \cite{Gro}, the direct sum $\oplus_{n=0}^\infty
K_{S_n}(X^n)\otimes \C$ of equivariant $K$-groups can be
identified with a Fock space of a Heisenberg algebra, cf.
\cite{Seg, Wa1}. In particular, the size of the cohomology group
$H^*(\Xn)$ coincides with that of the equivariant $K$-group
$K_{S_n}(\Xn)$ or the orbifold cohomology group of $X^n/S_n$.
(Equivariant $K$-groups are related to orbifold cohomology groups
by a decomposition theorem due to \cite{BC} and independently
\cite{Kuhn}). However the ``equivalence" on the level of ring
structures is more subtle.

In \cite{CR}, Chen and Ruan introduced an orbifold cohomology {\em
ring} for any orbifold. When $X$ is a projective surface with a
numerically trivial canonical class, Ruan \cite{Ru1} conjectured
that the orbifold cohomology ring of the symmetric product
$X^n/S_n$ is isomorphic to the cohomology ring of Hilbert scheme
$\Xn$. This has been established with some sign modification in
\cite{LS2, FG}. More precisely, a certain graded Frobenius algebra
$A^{[n]}$ was constructed in Lehn-Sorger \cite{LS2} based on a
graded Frobenius algebra $A$. The cohomology ring $H^*(\Xn)$ is
then shown to be isomorphic to $A^{[n]}$ for $A=H^*(X)$, which is
subsequently identified in \cite{FG} with the rational orbifold
cohomology ring $\orbsym$ with some modifications of signs in the
orbifold cup product, also see \cite{Uri}.

The proof of the ring isomorphism in \cite{LS2} is very ingenious
however quite indirect. It made use of earlier results on Hilbert
schemes as well as an observation in \cite{FW} on the relation
between Goulden's operator on the symmetric groups and Lehn's
operator in Hilbert schemes. In particular, much remains to be
understood about the finer structures of the ring $\orbsym$ on its
own for a general complex manifold $X$ or about any direct
connections between $\orbsym$ and vertex algebras.

The goal of the present paper is to develop systematically a
(vertex) algebraic approach to study the ring $\orbsym$ for a
closed complex manifold $X$ of even dimension in a self-contained
manner which is {\em independent} of (but parallel to) the
development on Hilbert schemes. By examining closely the arguments
in [{\bf LQW1-4}], we observe that the results obtained therein
follow in an axiomatic way from several key constructions and
statements obtained in \cite{Na2, Lehn, LQW1}, no matter how
difficult or how long the proofs could be. (This does not mean
that we could formulate and prove these results in Hilbert schemes
in any easier way.) Thus we formulate several axioms to formalize
the setup in order to make it applicable in different situations.
We obtain the corresponding key constructions and prove the
corresponding statements in $\orbsym$. These enable us to apply
the axioms to obtain the counterparts in the setup of symmetric
products of the main results in [{\bf LQW1-4}] for Hilbert
schemes.

Let us explain in more detail. We first note that a Heisenberg
algebra acting on $\Fock = \oplus_{n=0}^\infty \orbsym$ is
available by easily reformulating the equivariant $K$-group
construction in \cite{Seg, Wa1}. We then introduce the cohomology
classes $O^k(\alpha,n) \in \orbsym$ using the $k$-th power sum of
Jucys-Murphy elements in the symmetric groups \cite{Juc, Mur}. We
define the operator $\mathfrak O^k(\alpha) \in \End(\Fock)$ to be
the orbifold cup product with $O^k(\alpha,n)$ in $\orbsym$ for
each $n$. (As we shall see, the operators $\mathfrak O^k(\alpha)$
turn out to be the counterpart of the Chern character operators in
\cite{LQW1}.) The idea of relating Jucys-Murphy elements to vertex
operators has been used in \cite{LT} in the study of class
functions of symmetric groups. In particular, when $X$ is a point,
the operator $\mathfrak O^1$ reduces to the Goulden's operator
\cite{Gou} which admits a vertex operator interpretation
\cite{FW}.

By studying the interaction of the operators $\mathfrak
O^k(\alpha)$ and Heisenberg operators, we are able to verify all
the axioms. Therefore, as formal consequences of the
axiomatization, we obtain the counterparts in the symmetric
product setup of all the main results for the cohomology rings of
Hilbert schemes in [{\bf LQW1-4}]. Namely, we prove that
$O^k(\alpha,n)$, as $0 \le k <n$ and $\alpha$ runs over a linear
basis of $H^*(X)$, form a set of ring generators of $\orbsym$. We
also show that there is another set of ring generators in terms of
Heisenberg algebra generators. We establish the stability of the
ring $\orbsym$, which tells us in what sense the orbifold cup
product on $\orbsym$ is independent of $n$. We further obtain a
description of the operators $\mathfrak O^k(\alpha)$ as the zero
mode of a certain explicit vertex operator when $X$ has a positive
dimension. (Such a description has been given in \cite{LT} when
$X$ is a point. There is however a remarkable difference between
these two cases.) The components of these vertex operators
generate a $\W$ algebra acting on $\Fock$. The description of the
operator $\mathfrak O^k(\alpha)$ as the zero mode of a certain
explicit vertex operator provides us a new way to construct a
sequence of Frobenius algebra $\mathfrak F_A^n$ starting from a
Frobenius algebra $A$, cf. Remark~\ref{orbrem:zeromode}. (Compare
with the different construction in \cite{LS2} of the Frobenius
algebra $A^{[n]}$).

For a global quotient $M/G$, we introduce a deformed orbifold cup
product on $H^*_{\text{orb}}(M/G)$ (actually first on $H^*(M,G)$)
depending on a rational (or complex, if we consider the orbifold
cohomology group with $\C$-coefficient) parameter $t$. This
reduces to the original construction in \cite{CR} for $t=1$, and
to the construction of \cite{FG} (also cf. \cite{LS2}) for $t=-1$.
In the case of symmetric products, we explain how the parameter
$t$ is reflected in terms of some modifications on Heisenberg
algebra and vertex operators. By comparing our results on the
symmetric products with results on Hilbert schemes in \cite{LQW1},
we obtain a new proof of the ring isomorphism between $\Hxn$ and
$\orbsym$ (with the sign modification) for $X$ with a numerically
trivial canonical class.

It turns out that the $t$-family of deformed ring structures on
$H^*_{\text{orb}}(M/G,\C)$ {\em with $\C$-coefficient} are
isomorphic for all nonzero $t$. As a consequence, using $\C$
instead of $\mathbb Q$ as coefficients for (orbifold) cohomology
groups, there exists a ring isomorphism between cohomology ring of
$\Xn$ and the original orbifold cohomology ring of $X^n/S_n$ for
$X$ with a numerically trivial canonical class. This supports the
original conjecture of Ruan \cite{Ru1} on the ring isomorphism of
hyperkahler resolutions, if we insist on using $\C$ rather than
$\mathbb Q$ as the cohomology coefficients.

As observed in \cite{Wa1}, for a given complex manifold $X$ with a
finite group $\G$ action, the product $X^n$ affords a natural
action of the wreath product $\Gn$ (which is a finite group given
by the semidirect product $\G^n \rtimes S_n$). Further, the
quotient $X^n/\Gn$ can be identified with the symmetric product of
the orbifold $X/\G$. The results on the orbifold cohomology ring
of symmetric products in this paper will be generalized elsewhere
to the symmetric products of a general orbifold.

We find it amazing to have such a wonderful dictionary between
results in Hilbert schemes and in symmetric products. To some
extent, the results on symmetric products are simpler since the
canonical class does not play a role here (besides, $X$ needs not
to be a surface). We present a partial dictionary in a table near
the end of the paper. It is also instructive to compare with
another dictionary table in \cite{Wa2} between Hilbert schemes and
wreath products.

Our results in this paper may shed light on the understanding of
the difference between the rings $\Hxn$ and $\orbsym$ (with/out
the sign changes) even when the projective surface has a
nontrivial canonical class. A special case of a conjecture made in
\cite{Ru2} says that there exists a ring isomorphism between
$\orbsym$ and $\Hxn$ with a quantum corrected product. A partial
verification has been made by the Gromov-Witten 1-point function
computation for $\Xn$ in \cite{LQ}. Our axiomatization provides a
possible strategy for checking Ruan's conjectural ring
isomorphism. Namely, we can use the quantum corrected product to
replace the usual cup product to introduce operators analogous to
the Chern character operators, and then try to understand their
interaction with the {\em usual} Heisenberg algebra, and then
compare with the results on the ring $\orbsym$ obtained in this
paper.

The paper is organized as follows. In Sect.~\ref{sec:hilb}, we
review the key constructions and statements in Hilbert schemes
which are responsible for the further results. We put an emphasis
on the axiomatic nature of these results. We also obtain a
variation of Lehn's theorem in relating the Chern class of certain
tautological bundles on $\Xn$ and Heisenberg generators, and point
out an interesting corollary. In Sect.~\ref{sec:orbkey}, we
present the corresponding key constructions and prove the
corresponding statements in symmetric products. The reader should
compare the constructions in Sect.~\ref{sec:orbkey} and in
Sect.~\ref{sec:hilb}. In Sect.~\ref{sec:orbformal}, we formulate
formal consequences in symmetric products of the results in the
previous section, which are the counterparts of earlier results on
Hilbert schemes. In Sect.~\ref{sec:sign}, we explain how a
deformation of the orbifold cup product $\orbsym$ is reflected in
terms of modified Heisenberg algebra and vertex operators. As a
corollary, we obtain a new proof of the modified Ruan's conjecture
on the ring isomorphism between $\Hxn$ and $\orbsym$ when $X$ has
a numerically trivial canonical class. Further, we show that
Ruan's original conjecture holds if we use $\C$ as cohomology
coefficient. In Sect.~\ref{sec:open}, we list several open
questions for further research.

{\bf Convention.} All the (orbifold) cohomology groups/rings are
assumed to have $\mathbb Q$-coefficients unless otherwise
specified.
\section{The cohomology ring of Hilbert schemes}
\label{sec:hilb}
\subsection{Hilbert schemes of points on surfaces}

Let $X$ be a smooth projective complex surface with the canonical
class $\KX$ and the Euler class $e$, and $\Xn$ be the Hilbert
scheme of points in $X$. We define a bilinear form
$$ (\alpha,\beta) =\int_X \alpha\beta,\qquad \alpha,\beta \in
H^*(X). $$
An element in $\Xn$ is represented by a length-$n$ $0$-dimensional
closed subscheme $\xi$ of $X$. For $\xi \in \Xn$, let $I_{\xi}$ be
the corresponding sheaf of ideals. It is well known that $\Xn$ is
smooth. Sending an element in $\Xn$ to its support in the
symmetric product $X^n/S_n$, we obtain the Hilbert-Chow morphism
$\pi_n: \Xn \rightarrow X^n/S_n$, which is a resolution of
singularities. Define the universal codimension-$2$ subscheme:
\begin{eqnarray*}
{ \mathcal Z}_n=\{(\xi, x) \subset \Xn\times X \, | \, x\in
 {\rm Supp}{(\xi)}\}\subset \Xn\times X.
\end{eqnarray*}
Denote by $p_1$ and $p_2$ the projections of $\Xn \times X$ to
$\Xn$ and $X$ respectively. Let
\begin{eqnarray*}
\Fockh = \bigoplus_{n=0}^\infty H^*(\Xn)
\end{eqnarray*}
be the direct sum of total cohomology groups (with $\mathbb
Q$-coefficient) of the Hilbert schemes $\Xn$.
\subsection{The Heisenberg Algebra}

Nakajima and Grojnowski \cite{Na1, Gro} constructed geometrically
a Heisenberg algebra which acts irreducibly on $\Fockh$ with
generators $\mathfrak a_n(\alpha)$, $n\in \Z, \alpha \in H^*(X)$.
Below we recall the construction of Nakajima \cite{Na2}.

For $m \ge 0$ and $n > 0$, let $Q^{[m,m]} = \emptyset$ and define
$Q^{[m+n,m]}$ to be the closed subset:
$$\{ (\xi, x, \eta) \in X^{[m+n]} \times X \times X^{[m]} \, | \,
\xi \supset \eta \text{ and } \mbox{Supp}(I_\eta/I_\xi) = \{ x \}
\}.$$

Let $n \ge 0$. The linear operator $\mathfrak a_{-n}(\alpha) \in
\End(\Fockh)$ with $\alpha \in H^*(X)$ is defined by
$$\mathfrak a_{-n}(\alpha)(a) = \tilde{p}_{1*}([Q^{[m+n,m]}] \cdot
\tilde{\rho}^*\alpha \cdot \tilde{p}_2^*a)$$
for $a \in H^*(X^{[m]})$, where $\tilde{p}_1, \tilde{\rho},
\tilde{p}_2$ are the projections of $X^{[m+n]} \times X \times
X^{[m]}$ to $X^{[m+n]}, X, X^{[m]}$ respectively. Define
$\mathfrak a_{n}(\alpha) \in \End(\Fockh)$ to be $(-1)^n$ times
the operator obtained from the definition of $\mathfrak
a_{-n}(\alpha)$ by switching the roles of $\tilde{p}_1$ and $
\tilde{p}_2$. We often refer to $\mathfrak a_{-n}(\alpha)$ (resp.
$\mathfrak a_n(\alpha)$) as the {\em creation} (resp. {\em
annihilation})~operator.

\begin{theorem} \label{hilb:heis}
 The operators $\mathfrak a_n(\alpha) \in \End(\Fockh)$
 $(n\in \Z, \alpha \in H^*(X))$
 generate a Heisenberg (super)algebra with commutation relations
 given by
\begin{eqnarray*}
  [ \mathfrak a_{m} (\alpha), \mathfrak a_{n} (\beta)]
  = -m \delta_{m,-n} (\alpha, \beta) \cdot \text{\rm Id}_{\Fockh}
\end{eqnarray*}
where $n,m\in \Z,\;\alpha, \beta \in H^*(X)$. Furthermore,
$\Fockh$ is an irreducible representation of the Heisenberg
algebra with the vacuum vector $\vac =1 \in H^*(pt) \cong \C$.
\end{theorem}
The commutator above is understood in the super sense according to
the parity of cohomology classes $\alpha, \beta$ involved.

We define the following cohomology class in $H^*(\Xn)$
\cite{LQW2}:
$$B_i(\g, n) = \frac{1}{(n-i-1)!}
\cdot \mathfrak a_{-i-1}(\gamma) \mathfrak
a_{-1}(1_X)^{n-i-1}\vac.$$
\subsection{The Chern character operator $\mathfrak G^k(\alpha)$}

In \cite{LQW1} we introduced the cohomology classes $G(\g,n)$,
$G_k(\g,n)$, and the operators $\mathfrak G_k(\alpha)$ etc (also
cf. \cite{Lehn}).

For $n \ge 0$ and a homogeneous class $\gamma \in H^*(X)$, let
$|\gamma| = s$ if $\gamma \in H^s(X)$, and let $G_i(\gamma, n)$ be
the homogeneous component in $H^{|\gamma|+2i}(\Xn)$ of
\[
 G(\gamma, n) = p_{1*}(\ch({\mathcal O}_{{\mathcal Z}_n}) \cdot
 p_2^*{\rm td}(X) \cdot p_2^*\gamma) \in H^*(\Xn)
\]
where $\ch({\mathcal O}_{{\mathcal Z}_n})$ denotes the Chern
character of the structure sheaf ${\mathcal O}_{{\mathcal Z}_n}$
and ${\rm td}(X) $ denotes the Todd class. Here and below we omit
the Poincar\'e duality used to switch a homology class to a
cohomology class and vice versa. We extend the definition of
$G_i(\gamma, n)$ to an arbitrary class $\gamma \in H^*(X)$ by
linearity. It turns out to be more convenient to introduce a
normalized class
$$G^k(\g,n) :=k! \cdot G_k(\g, n).$$
It was proved in \cite{LQW1} that the cohomology ring of $\Xn$ is
generated by the classes $G^{i}(\gamma, n)$ where $0 \le i < n$
and $\gamma$ runs over a linear basis of $H^*(X)$.

The {\it Chern character operator} ${\mathfrak G}^i(\gamma) \in
\End({\Fockh})$ is defined to be the operator acting on the
component $H^*(\Xn)$ by the cup product with $G^i(\gamma, n)$ for
every $n \ge 0$. We introduce a formal variable $\hbar$ (here and
in other places later on) and let $\Gh (\g) = \sum_{i \ge 0}
\frac{\hbar^i}{i!} \cdot {\mathfrak G}^i(\g).$ A convenient way is
to regard $\hbar$ as having `cohomology degree' $-2$ so $\Gh (\g)$
becomes homogeneous of degree $|\g|$.

Let $\mathfrak d = \mathfrak G^1(1_X)$ where $1_X$ is the
fundamental cohomology class of $X$. The operator $\mathfrak d$
was first introduced in \cite{Lehn} and plays an important role in
the theory. For a linear operator $\mathfrak f \in \End(\Fockh)$,
define its {\it derivative} $\mathfrak f'$ by $\mathfrak f' =
[\mathfrak d, \mathfrak f]$. The higher derivative $\mathfrak
f^{(k)}$ is defined inductively by $\mathfrak f^{(k)} = [\mathfrak
d, \mathfrak f^{(k-1)}]$.
\subsection{Vertex operators}

We define the normally ordered product $:\mathfrak
a_{m_1}\mathfrak a_{m_2}:$ to be $\mathfrak a_{m_1}\mathfrak
a_{m_2}$ when $m_1 \le m_2$ and $\mathfrak a_{m_2}\mathfrak
a_{m_1}$ when $m_1 > m_2$. We denote
$$\tau_{k*}: H^*(X) \to
H^*(X^k) \cong H^*(X)^{\otimes k},$$ where $k \ge 1$, is the
linear map induced by the diagonal embedding $\tau_k: X \to X^k$,
and $\mathfrak a_{m_1} \cdots \mathfrak
a_{m_k}(\tau_{k*}(\alpha))$ denotes $\sum_j \mathfrak
a_{m_1}(\alpha_{j,1}) \cdots \mathfrak a_{m_k}(\alpha_{j,k})$ if
we write $\tau_{k*}\alpha = \sum_j \alpha_{j,1} \otimes \cdots
\otimes \alpha_{j, k}$ via the K\"unneth decomposition of
$H^*(X^k)$. We will simply write $\tau_*$ for $\tau_{k*}$ when
there is no confusion.

Our convention of vertex operators or fields is to write them in a
form
\begin{eqnarray*}
\phi(z) =\sum_n \phi_n z^{-n-\Delta}
\end{eqnarray*}
where $\Delta$ is the conformal weight of the field $\phi (z)$.
Define the {\it derivative field}
\begin{eqnarray*}
\partial \phi(z) =\sum_n (-n-\Delta) \phi_n z^{-n-\Delta-1}.
\end{eqnarray*}
We define the normally ordered product $:\phi_1(z) \cdots
\phi_k(z):$ as usual (cf. e.g. \cite{Kac}).

For $\alpha \in H^*(X)$, we define a vertex operator
$a(\alpha)(z)$ by putting
\begin{eqnarray*}
a(\alpha)(z) = \sum_{n \in \Z} \mathfrak a_{n}(\alpha) z^{-n-1}.
\end{eqnarray*}
The field $:a(z)^p:(\tau_{*}\alpha)$ is defined to be $\sum_i
:a(\alpha_{i,1})(z) a(\alpha_{i,2})(z) \cdots a(\alpha_{i,p})(z):$
if we write $\tau_{p*} \alpha =\sum_i \alpha_{i,1}\otimes
\alpha_{i,1}\otimes \ldots \otimes \alpha_{i,p} \in
H^*(X)^{\otimes p}$. We rewrite $:a(z)^{p}:(\ta)$ componentwise as
\begin{eqnarray*}
:a(z)^{p}: (\ta)= \sum_m :a^{p}:_m (\ta)\; z^{-m-p},
\end{eqnarray*}
where $:a^{p}:_m (\ta) \in \End(\Fockh)$ is the coefficient of
$z^{-m-p}$ (i.e. the $m$-th Fourier component of the field
$:a(z)^{p}: (\ta)$), and maps $H^*(X^{[n]})$ to $H^*(X^{[n+m]})$.
Similarly, for $r \ge 1$, we can define the field $:(\partial^r
a(z)) a(z)^{p-1}:(\tau_{*}\alpha)$, and define the operator
$:(\partial^r a)a^{p-1}:_m(\tau_{*}\alpha)$ as the coefficient of
$z^{-m- r-p}$ in $:(\partial^r a(z))a(z)^{p-1}:(\tau_{*}\alpha)$.
\subsection{Interactions between Heisenberg algebra and $\mathfrak G(\alpha)$}

The following theorem is a variation of Lemma~5.8, \cite{LQW1},
which generalizes Theorem~4.2, \cite{Lehn}. Clearly the two
identities in the following theorem are equivalent.

\begin{theorem}  \label{hilb:comm}
 Let $\g, \alpha \in H^*(X)$. Then we have
 \begin{eqnarray*}
[\Gh(\g), \mathfrak a_{-1}(\alpha)]
 &=& \exp(\hbar \cdot {\rm ad}\;{\mathfrak d}) (\mathfrak a_{-1}(\g \alpha)) \\
 {[} \mathfrak G^k(\g), \mathfrak a_{-1}(\alpha)]
 &=& \mathfrak a^{(k)}_{-1}(\g\alpha),\quad k \ge 0.
\end{eqnarray*}
\end{theorem}

\begin{theorem}  \label{hilb:cubic}
For $\alpha \in H^*(X)$, we have
 \begin{eqnarray}  \label{hilbeq:general}
 \mathfrak G^1(\alpha) =-\frac16 :\mathfrak a^3:_0
(\tau_{*}\alpha)
 -\sum_{n>0} \frac{n-1}2 :\mathfrak a_n \mathfrak a_{-n}:
 (\tau_*(K\alpha)).
 \end{eqnarray}
 In particular, for a surface $X$ with numerically trivial
canonical class, we have
\begin{eqnarray}  \label{hilbeq:cubic}
 {\mathfrak d} =-\frac16 :\mathfrak a^3:_0 (\tau_{*}1_X).
\end{eqnarray}
\end{theorem}

\begin{proof}
Observe that both sides of (\ref{hilbeq:general}) annihilates the
vacuum vector $\vac$. To prove (\ref{hilbeq:general}), it suffices
to show that the commutators of both sides of
(\ref{hilbeq:general}) with the operators $\mathfrak a_n(\beta)$,
$n\in\Z, \beta \in H^*(X)$, coincide. It was shown in \cite{LQW1}
that (i.e. the {\em transfer property})
 $$[\mathfrak G^1(\alpha),\mathfrak a_n(\beta)]
=[\mathfrak G^1(1_X),\mathfrak a_n(\alpha\beta)] =\mathfrak
a_n'(\alpha\beta),$$ while $\mathfrak a_n'(\alpha\beta)$ was
computed in \cite{Lehn}, Theorem 3.10. We can easily compute the
commutator of the right hand side of (\ref{hilbeq:general}) with
$\mathfrak a_n(\beta)$ by using Lemma~3.1 in \cite{LQW3}. These
two commutators coincide.
\end{proof}
\subsection{Axiomatization}  \label{subsec:axiom}

We claim that Theorem~\ref{hilb:heis}, Theorem~\ref{hilb:comm},
and Theorem~\ref{hilb:cubic} encode all the information about the
ring structure of the cohomology ring $H^*(\Xn)$ for every $n$. In
fact, if we examine closely the proofs of all the main results in
[{\bf LQW1-4}], we see that these three statements (or sometimes
some weaker form of Theorem~\ref{hilb:cubic}) were used
effectively, together with numerous standard properties concerning
$\tau_{k*}$, the Heisenberg generators etc. The results therein
include the theorems on cohomology ring generators in \cite{LQW1,
LQW2}, the universality and stability theorems in \cite{LQW3}, and
the connection with $\W$ algebras in \cite{LQW4}.

This observation leads to a strategy which allows us to treat
similar cases in an axiomatic manner. Assume that ({\bf A1}) there
exists a sequence of (finite-dimensional) graded Frobenius
algebras $A^{[n]}$ ($n \ge 0$) such that $A=A^{[1]}$. ({\bf A2})
the direct sum $\oplus_n A^{[n]}$ affords the structure of a Fock
space of a Heisenberg algebra. ({\bf A3}) There exists a sequence
of elements $G^k(\alpha,n) \in A^{[n]}$ depending on $\alpha \in
A$ (linearly) and a non-negative integer $k$, which can be used to
define operators $\mathfrak G^k(\alpha)$. The operators $\mathfrak
G^k(\alpha)$, $\mathfrak G^1(1_A)$ and the Heisenberg generators
satisfy the relations as in Theorem~\ref{hilb:comm} and
(\ref{hilbeq:cubic}). (Note here that in general there is no
counterpart of the $K$-term in Theorem~\ref{hilb:cubic}, so we
have put a stricter axiom here. Variations of this are allowed in
different setup.)

The axioms ({\bf A1-A3}) suffice us to use the same approach as in
[{\bf LQW1-4}] to obtain similar theorems. (For results in
\cite{LQW3,LQW4} the fact $e^2=0$ is used. This requires a
separate treatment when $A\cong \C$ which corresponds to the
cohomology ring of $X =pt$.) In fact, many formulas in the Hilbert
scheme side are to be simplified by discarding all the terms which
involve in a manifest way with the canonical class $K$.

We will follow this axiomatic route when we treat the orbifold
cohomology ring of the symmetric products in the sections below.
\subsection{Tautological bundles and Heisenberg generators}

Given a line bundle $L$ over $X$, we obtain a rank $n$ vector
bundle $L^{[n]}=(p_1|_{\mathcal Z_n})_* ( p_2^*L|_{\mathcal Z_n})$
over $\Xn$, where we recall that $p_1$ and $p_2$ are the
projections from $\Xn \times X$ to $\Xn$ and $X$ respectively, and
$\mathcal Z_n$ is the universal subscheme of $\Xn \times X$. We
introduce the generating function for Chern classes:
$${c}_{\hbar}(L^{[n]}) =\sum_{i \ge 0} {c}_i(L^{[n]}) \hbar^i.$$

The following theorem is a variation of Theorem~4.6 in
\cite{Lehn}. We can prove it in a way parallel to the Variant 2 of
the proof of Lehn's theorem, once we make sure to put a suitable
power of $\hbar$ in the right places.

\begin{theorem}  \label{th:generating}
Let $L$ be a line bundle on $X$. Then
\begin{eqnarray} \label{eq_general}
 \sum_{n=0}^\infty c_\hbar (L^{[n]}) z^n
 =\exp \left( \sum_{r\ge 1}\frac{(-\hbar)^{r-1}}r
   \mathfrak a_{-r}(c_\hbar(L)) z^r \right) \cdot\vac.
\end{eqnarray}
\end{theorem}

When setting $\hbar=1$, we recover Theorem~4.6, \cite{Lehn}:
 $$\sum_{n=0}^\infty c(L^{[n]}) z^n
 =\exp\left( \sum_{ r\ge 1} \frac{(-1)^{r-1}}r \mathfrak a_{-r}(c(L)) z^r\right)\cdot
 \vac.$$

Besides making various gradings involved more transparent, one
advantage of our Theorem~\ref{th:generating} to Theorem~4.6 in
\cite{Lehn} is that we may set the value of $\hbar$ to be a number
other than 1. Denote by $L^{[n]\vee}$ the dual bundle of
$L^{[n]}$. The Chern classes of a vector bundle and its dual are
related to each other: $c_k(L^{[n] \vee}) =(-1)^k c_k(L^{[n]})$,
$k \ge 0$. Setting $\hbar=-1$ in Eq.~(\ref{eq_general}), we arrive
at the following new observation.

\begin{corollary}
 Let $L$ be a line bundle on $X$. Then
 $$\sum_{n=0}^\infty c(L^{[n]\vee}) z^n
 =\exp\left( \sum_{ r\ge 1} \frac1r \mathfrak a_{-r}(c(L^\vee)) z^r\right)\cdot
 \vac.$$
\end{corollary}
\section{The orbifold cohomology ring of symmetric products I}
\label{sec:orbkey}
\subsection{Generalities on orbifold cohomology rings}

Let $M$ be a complex manifold of complex dimension $d$ with a
finite group $G$ action. Following \cite{BBM, BC, Kuhn}, we
introduce the space

$$M\diamond G =\{(g, x) \in G \times M\mid gx=x\} =\bigsqcup_{g \in G} M^g,$$
and $G$ acts on $M\diamond G$ naturally by $h.(g,x) =(hgh^{-1},
hx)$. As a vector space, we define $H^*(M,G)$ to be the cohomology
group of $M\diamond G$ with rational coefficient (cf. \cite{FG}),
or equivalently,

$$H^*(M,G) = \bigoplus_{g\in G} H^*(M^g).$$
The space $H^*(M,G)$ has a natural induced $G$ action, which is
denoted by $\ad \, h : H^*(M^g) \rightarrow H^*(M^{hgh^{-1}})$. As
a vector space, the orbifold cohomology group
$H^*_{\text{orb}}(M/G)$  is the $G$-invariant part of $H^*(M,G)$,
which is isomorphic to
 $$\bigoplus_{[g] \in G_*} H^*(M^g /Z(g))$$
where $G_*$ denotes the set of conjugacy classes of $G$ and
$Z(g)=Z_G(g)$ denotes the centralizer of $g$ in $G$.

For $g\in G$ and $x \in M^g$, write the eigenvalues of the action
of $g$ on the complex tangent space $TM_x$ to be $\mu_k =e^{2\pi i
r_k}$, where $0 \le r_k <1$. The {\em degree shift number} (or
{\em age}) is the rational number $F^g_x =\sum_{k=1}^d r_k$, cf.
\cite{Zas}. It depends only on the connected component $Z$ which
contains $x$, so we can denote it by $F^g_Z$. Then associated to a
cohomology class in $H^r(Z)$, we assign the corresponding element
in $H^*(M,G)$ (and thus in $H^*_{\text{orb}}(M/G)$) the degree
$r+2 F^g_Z$.

A ring structure on $H^*_{\text{orb}}(M/G)$ was introduced by Chen
and Ruan \cite{CR}. This was subsequently clarified in \cite{FG}
by introducing a ring structure on $H^*(M,G)$ first and then
passing to $H^*_{\text{orb}}(M/G)$ by restriction. We shall use
$\circ$ to denote this product. The ring structure on $H^*(M,G)$
is degree-preserving, and has the property: $\alpha \circ \beta $
lies in $H^*(M^{gh})$ for $\alpha \in H^*(M^g)$ and $\beta \in
H^*(M^h)$.

For $1\in G$, $H^*(M^1/Z(1)) \cong H^*(M/G)$, and thus we can
regard $\alpha \in H^*(M/G)$ to be $\alpha \in
H^*_{\text{orb}}(M/G)$ by this isomorphism. Also given $a
=\sum_{g\in G} a_g g$ in $\mathbb Q [G]$ (resp. $\mathbb Q
[G]^G$), we may regard $a$ as an element in $H^*(M,G)$ (resp.
$H^*_{\text{orb}}(M/G)$) whose component in each $H^*(M^g)$ is
$a_g \cdot 1_{M^g} \in H^0(M^g)$.

If $K$ is a subgroup of $G$, then we can define the restriction
map from $ H^*(M,G)$ to $H^*(M,K)$ by projection to the component
$\oplus_{g \in K} H^*(M^g)$ which, when restricted to the
$G$-invariant part, induces naturally a {\em degree-preserving}
linear map $ \Res^G_K: H^*_{\text{orb}}(M/G) \rightarrow
H^*_{\text{orb}}(M/K)$. We define the induction map

 $$\Ind_K^G:H^*(M,K) \rightarrow H^*_{\text{orb}}(M/G)$$
by sending $\alpha \in H^*(M^h),$ where $h \in K$, to
$$\Ind_K^G(\alpha) = \frac1{|K|}\sum_{g \in G}  \ad g
(\alpha ) .$$
Note that $\Ind_K^G (\alpha)$ is clearly $G$-invariant. When
restricted to the invariant part, we obtain a {\em
degree-preserving} linear map $ \Ind_K^G: H^*_{\text{orb}}(M/K)
\rightarrow H^*_{\text{orb}}(M/G)$. We often write the restriction
and induction maps as $\Res_K, \Res$ and $ \Ind^G, \Ind$, when the
groups involved are clear from the context. In particular, when
$M$ is a point, $H^*_{\text{orb}}(pt/G)$ reduces to the
Grothendieck ring $R_{\mathbb Q}(G)$ of $G$, and we recover the
induction and restriction functors in the theory of finite groups.
\subsection{The Heisenberg algebra}

Let $X$ be a closed complex manifold of complex dimension $d$. Our
main objects are the orbifold cohomology ring $\orbsym$, and the
non-commutative ring $H^*(X^n,S_n)$. We denote
 $$\Fock =\bigoplus_{n=0}^\infty \orbsym. $$

We introduce a linear map $\omega_n: H^*(X) \rightarrow
H^*_{\text{orb}}(X^n/S_n)$ as follows: given $\alpha \in H^r(X)$,
we denote by $\omega_n(\alpha) \in
H^{r+d(n-1)}_{\text{orb}}(X^n/S_n)$ the element associated to
$n\alpha$ by the isomorphism $H^*((X^n)^{\sigma_n }) \cong H^*(X)$
for any permutation $\sigma_n$ in the conjugacy class $[n]\in
(S_n)_*$ which consists of the $n$-cycles. We also define $\ch_n:
H^*_{\text{orb}}(X^n/S_n) \rightarrow H^*(X)$ as the composition
of the isomorphism $H^*((X^n)^{\sigma_n}) \cong H^*(X)$ with the
projection from $H^*_{\text{orb}}(X^n/S_n)$ to
$H^*((X^n)^{\sigma_n})$.

Let $\alpha \in H^*(X)$. For $n>0$, we define the creation operator
$\mathfrak p_{-n}(\alpha)\in \End (\Fock)$ given by
the composition ($k \ge 0$):

\begin{eqnarray*}
 H^*_{\text{orb}}(X^k/S_k)
 &\stackrel{\omega_n(\alpha)\otimes \cdot}{\longrightarrow}&
 H^*_{\text{orb}}(X^n/S_n) \bigotimes H^*_{\text{orb}}(X^k/S_k)
 \\
 &\stackrel{\cong}{\longrightarrow}&
 H^*_{\text{orb}}(X^{n+k}/(S_n\times S_k))
 \stackrel{\Ind}{\longrightarrow}
 H^*_{\text{orb}}(X^{n+k}/S_{n+k}) ,
\end{eqnarray*}
and the annihilation operator $\mathfrak p_{n}(\alpha)\in \End
(\Fock)$ given by the composition ($k \ge 0$):

\begin{eqnarray*}
 H^*_{\text{orb}}(X^{n+k}/S_{n+k})
 & \stackrel{\text{Res}}{\longrightarrow} &
 H^*_{\text{orb}}(X^{n+k}/(S_n \times S_{k}))  \\
 & \stackrel{\cong}{\longrightarrow} &
 H^*_{\text{orb}}(X^n/S_n)
 \bigotimes  H^*_{\text{orb}}(X^{k}/ S_{k}) \\
 &\stackrel{\ch_n}{\longrightarrow} &
 H^*(X) \bigotimes H^*_{\text{orb}}(X^{k}/S_{k})
 \stackrel{(\alpha, \cdot)}{\longrightarrow}
 H^*_{\text{orb}}(X^{k}/S_{k}).
\end{eqnarray*}
We also set $\mathfrak p_0(\alpha) =0$.

\begin{theorem} \label{orb:heis}
 The operators $\mathfrak p_n(\alpha) \in \End(\Fock)$
$(n\in \Z, \alpha \in H^*(X))$
 generate a Heisenberg (super)algebra with commutation relations
 given by
\begin{eqnarray*}
  [ \mathfrak p_{m} (\alpha), \mathfrak p_{n} (\beta)]
  = m \delta_{m,-n} (\alpha, \beta) \cdot \text{\rm Id}_{\Fock}
\end{eqnarray*}
where $n,m \in \Z, \;\alpha,\beta \in H^*(X)$. Furthermore,
$\Fock$ is an irreducible representation of the Heisenberg algebra
with the vacuum vector $\vac =1 \in H^*(pt) \cong \C$.
\end{theorem}
This theorem can be proved in the same way as an analogous theorem
formulated by using the equivariant $K$-group $K_{S_n}(X^n)\otimes
\C$. This analogous theorem was established in \cite{Seg} (see
\cite{Wa1}, Theorem~4 and its proof for detail). In general,
equivariant $K$-groups are related to orbifold cohomology groups
by a decomposition theorem \cite{Kuhn, BC}. Note that there is a
(fundamental!) sign difference in the two commutators of
Theorems~\ref{hilb:heis} and \ref{orb:heis}.

In particular, for a given $y \in
H^{|y|}_{\text{orb}}(X^{n-1}/S_{n-1})$, by the definition of
$\mathfrak p_{-1}(\alpha)$ (where $\alpha \in H^{|\alpha|}(X)$)
and the induction map, we can write that

$$\mathfrak p_{-1}(\alpha) (y) = \frac1{(n-1)!} \sum_{g \in
S_n} \ad g\, (\alpha \otimes y) = \frac{(-1)^{|\alpha|\cdot
|y|}}{(n-1)!} \sum_{g \in S_n} \ad g \,(y \otimes \alpha).$$

For $0 \le i <n$, we introduce the following cohomology class in
$H^*(\Xn)$:
$$P_i(\g, n) = {1 \over (n-i-1)!}
\cdot \mathfrak p_{-i-1}(\gamma) \mathfrak
p_{-1}(1_X)^{n-i-1}\vac.$$
\subsection{Jucys-Murphy elements}

For a permutation $\sigma \in S_n$ of cycle type given by a
partition $\lambda =(\lambda_1, \lambda_2, \ldots, \lambda_l)$ of
length $\ell =\ell(\lambda)$, we denote $d(\sigma) =d(\la) =n
-l(\lambda)$. Let $\ch: \oplus_{n=0}^\infty R(S_n) \rightarrow
\Lambda$ be the Frobenius characteristic map from the direct sum
of (complex) class functions on the symmetric group $S_n$ to the
ring $\Lambda$ of symmetric functions in infinitely many
variables, cf. \cite{Mac}. Denote by $\eta_n$ and $\varepsilon_n$
the trivial and alternating characters of $S_n$. Then $\ch$ sends
$\eta_n$ and $\varepsilon_n$ to the $n$-th complete and elementary
symmetric functions in $\Lambda$ respectively. We denote by $p_r$
the $r$-th power sum symmetric function.

Recall \cite{Juc, Mur} that the Jucys-Murphy elements $\xi_{j;n}$
of the symmetric group $S_n$ are defined to be the sums of
transpositions:
\begin{eqnarray*}
\xi_{j;n} = \sum_{i<j} (i,j), \quad j =1, \ldots, n.
\end{eqnarray*}
When it is clear from the text, we may simply write $\xi_{j;n}$ as
$\xi_j$. Denote by $\Xi_n$ the set $\{\xi_1, \ldots, \xi_n\}$.
According to Jucys, the $k$-th elementary symmetric function
$e_k(\Xi_n)$ of $\Xi_n = \{\xi_1, \ldots, \xi_n\}$ is equal to the
sum of all permutations in $S_n$ having exactly $(n-k)$ cycles.
Therefore, we obtain
\begin{eqnarray} \label{Jucys}
 \varepsilon_n
 = \sum_{\sigma \in S_n} (-1)^{d(\sigma)} \sigma
 = \sum_{k=0}^n (-1)^k e_k(\Xi_n)
 = \prod_{i=1}^n (1-\xi_i).
\end{eqnarray}
Denote by $\varepsilon_n(\hbar) =\prod_{i=1}^n (1-\hbar \xi_i)$,
where $\hbar$ is a formal parameter here and below. We have
\begin{eqnarray} \label{identity}
 \sum_{n=0}^\infty \ch (\varepsilon_n(\hbar)) z^n
 =\exp \left( \sum_{r\ge 1} (-\hbar)^{r-1}\frac{p_r}r z^r \right) .
\end{eqnarray}
Noting that $\varepsilon_n(1)$ (resp. $\varepsilon_n(-1)$)
coincides with the alternating character $\varepsilon_n$ (resp.
the trivial character $\eta_n$), we obtain two classical
identities involving $\eta_n, \varepsilon_n,$ and $p_r$ by setting
$\hbar =\pm 1$ in (\ref{identity}).
\subsection{The cohomology classes $\eta_n(\g)$ and $O^k(\alpha,n)$}

In the rest of this paper, we will assume that $X$ is a closed
complex manifold of {\em even} complex dimension $d$. Given $\g
\in H^*(X)$, we denote
$$\g^{(i)}= 1^{\otimes i-1}
\otimes \g \otimes 1^{\otimes n-i} \in H^*(X^n),$$
and regard it to be a cohomology class in $H^*(X^n,S_n)$
associated to the identity conjugacy class. We define $ \xi_i(\g)
:=\xi_i + \g^{(i)} \in H^*(X^n,S_n)$. We sometimes write
$\xi_i(\g)$ as $\xi_{i;n}(\g)$ to specify its dependence on $n$
when necessary.

\begin{definition}
Given $\g \in H^*(X)$, we define $\eta_n(\gamma)$ (resp.
$\varepsilon_n(\g)$) to be the cohomology class in $\orbsym$ whose
component associated to an element $\sigma$ in the conjugacy class
of partition $\lambda$ of $n$ is given by $\gamma^{\otimes
\ell(\lambda)} \in H^*((X^n)^\sigma)\cong H^{\otimes \ell
(\lambda)}$ (resp. by $(-1)^{d(\lambda)} \gamma^{\otimes
\ell(\lambda)}$). We further define an operator ${\bf \eta}(\g)$
(resp. ${\bf \varepsilon}(\g)$) in $\End (\Fock)$ by letting it
act on $\orbsym$ by the orbifold product with $\eta_n(\g)$ (resp.
$\varepsilon_n(\g)$) for every $n$.
\end{definition}

This definition is motivated by its counterpart in terms of
equivariant $K$-groups \cite{Seg, Wa1}. We can show that
$$\sum_{n=0}^\infty \eta_n(\g ) z^n
 =\exp\left( \sum_{ r\ge 1} \frac1r \mathfrak p_{-r}(\gamma) z^r\right)\cdot
 \vac.$$
Compare with \cite{Wa1}, Proposition~4.

\begin{proposition} \label{orbprop:eta}
Given $\g \in H^*(X)$, the orbifold cup product of the $n$
elements $\xi_i(\g)$ $(i=1,\ldots,n)$ in $H^*(X^n,S_n)$ lies in
$\orbsym$, and furthermore the following identity holds:

\begin{eqnarray} \label{eq:ident}
\eta_n(\g) =\prod_{i=1}^n \xi_i(\g) =  \xi_1(\g)\circ \xi_2(\g)
\circ \ldots \circ \xi_n(\g).
\end{eqnarray}
\end{proposition}

\begin{proof}
It suffices to prove (\ref{eq:ident}), since the first claim
follows from (\ref{eq:ident}) and the fact that $\eta_n(\g)$ is
$S_n$-invariant.

A typical monomial on the right-hand side of (\ref{eq:ident}) is
of the form
 $$(\xi_{i_1} \cdots \xi_{i_k}) \circ
(\gamma^{(j_1)}\cdots \gamma^{(j_{n-k})}) $$
where $i_1<\ldots <i_k$, $j_1 <\ldots <j_{n-k}$, and
$\{i_1,\ldots, i_k, j_1, \ldots, j_{n-k}\} =\{1, \ldots, n\}$.
Here we have used the observation that $\xi_{i_1} \circ \cdots
\circ \xi_{i_k}$ is just the usual multiplication of permutations
$ \xi_{i_1} \cdots \xi_{i_k}$ and $\gamma^{(j_1)} \circ\cdots
\circ \gamma^{(j_{n-k})}$ is just the ordinary cup product $
\gamma^{(j_1)}\cdots \gamma^{(j_{n-k})}$ in $H^*(X^n) \cong
H^*(X)^{\otimes n}$. Note that every cycle of each permutation
$\sigma$ appearing in $\xi_{i_1} \cdots \xi_{i_k}$ has exactly one
number which does not belong to ${i_1}, \ldots, {i_k}$, and in
addition, $\ell(\sigma) =n-k$ and $d(\sigma) =k$. Using the
definition of the orbifold cup product, we see that the product
$\sigma \circ \gamma^{(j_1)}\cdots \gamma^{(j_{n-k})}$ does not
involve the obstruction bundles (or the group defects are trivial
in the terminology of Lehn-Sorger) and equals $\g^{\otimes
\ell(\sigma)} \in H^{\otimes \ell(\sigma)} \cong H^*(
(X^n)^\sigma)$. This proves (\ref{eq:ident}).
\end{proof}

If we denote $ \varepsilon_n (\g,\hbar) = \prod_{i=1}^n (\g^{(i)}
-\hbar \xi_i),$ we have
\begin{eqnarray}  \label{eq:expon}
 \sum_{n=0}^\infty \varepsilon_n(\g,\hbar) z^n
 =\exp \left( \sum_{r\ge 1} (-\hbar)^{r-1}\frac{\mathfrak p_{-r}(\g)}r z^r \right) .
\end{eqnarray}

Regarding $\xi_i =\xi_i(0) \in H^*(X^n,S_n)$, we denote
$\xi_i^{\circ k} = \overbrace{ \xi_i \circ \ldots \circ\xi_i}^{k
~\rm{times}} \in H^*(X^n,S_n)$, and define $e^{-\xi_i} =\sum_{k \ge 0}
{1 \over k!} (-\xi_i)^{\circ k} \in H^*(X^n,S_n).$

\begin{definition}
For homogeneous $\alpha \in H^{|\alpha|}(X)$, we define the class
$O^k(\alpha,n) \in H^*_{\text{orb}}(X^n/S_n)$ to be
 $$O^k(\alpha,n) =\sum_{i=1}^n  (-\xi_i)^{\circ k} \circ \alpha^{(i)} \in
H^{dk+|\alpha|}_{\text{orb}}(X^n/S_n), $$
and extends linearly to all $\alpha \in H^*(X)$. We put
$O(\alpha,n) = \sum_{k\ge 0} \frac1{k!}O^k(\alpha,n) =\sum_{i=1}^n
e^{-\xi_i}\circ \alpha^{(i)},$ and put $O_{\hbar}(\alpha,n) =
\sum_{k\ge 0} \frac{\hbar^k}{k!}O^k(\alpha,n)$. We further define
the operator $\mathfrak O^k(\alpha) \in \End (\Fock)$ (resp.
$\mathfrak O (\alpha)$, or $\Oh(\alpha)$) to be the orbifold cup
product with $O^k(\alpha,n)$ (resp. $O(\alpha,n)$, or
$O_{\hbar}(\alpha,n)$) in $H^*_{\text{orb}}(X^n/S_n)$ for every $n
\ge 0$.
\end{definition}

\begin{remark}
We can see that $O^k(\alpha,n) \in H^*(X^n,S_n)$ is
$S_n$-invariant (and thus lies in $\orbsym$) as follows. For $\g
\in H^*(X)$, note that $e_j (\xi_1(\g), \ldots, \xi_n(\g))$
lies in $\orbsym$, where $e_j (\xi_1(\g), \ldots, \xi_n(\g))$
$(1\le j \le n)$ is the $j$-th elementary symmetric function in
$\xi_i(\g)$'s.
So $\orbsym$ contains all symmetric functions in $\xi_i(\g)$'s. In
particular, $O(e^{-\g}, n)  =\sum_i (e^{-\xi_i} \circ
(e^{-\g})^{(i)}) =\sum_i e^{-\xi_i(\g)}  \in \orbsym.$ Letting
$\g$ vary, we see that $O(\alpha,n)$ and similarly $O^k(\alpha,n)$
lie in $\orbsym.$
\end{remark}

In particular, $\mathfrak O^1(1_X) \in \End (\Fock)$ is the
generalized Goulden's operator \cite{Gou, FW}, which will be
denoted by $\mathfrak b$. The reason for our convention of putting
the sign in front of $\xi_i$ is to make the comparison with the
Hilbert scheme side easier. Also, the power sums of Jucys-Murphy
elements have been studied in \cite{LT} which corresponds to our
case when $X$ is a point.

Note that the generalized Goulden's operator $\mathfrak b$ is
defined to be the orbifold cup product with $O^1(1_X,n)$ in
$\orbsym$. Given an operator $\mathfrak f \in \End (\Fock)$, we
denote by $\mathfrak f ' = [\mathfrak b, f],$ and $\mathfrak
f^{(k+1)} = (\mathfrak f^{(k)})'.$ We have the following.

\begin{theorem}  \label{orb:cubic}
 We have ${\mathfrak b} =-\frac16 :\mathfrak p^3:_0
 (\tau_{*}1_X).$
\end{theorem}

\begin{remark} \label{orb:proofcubic}
This proposition is a counterpart of (\ref{hilbeq:cubic}). The
proof is essentially the same as the proof in the case when $X$ is
a point \cite{Gou, FW} (also cf. \cite{LS2}). For example, if we
look at the proof of Proposition~4.4, \cite{LS2}, the $\Delta_*$
and $e$ there should be replaced by our $\tau_*$ (which equals
$-\Delta_*$) and $-e$ respectively, since we are using the
orbifold cup product of \cite{CR}. Also compare
Proposition~\ref{orb:deformcubic} and
Remark~\ref{orb:deformproofcubic} below.
\end{remark}
\subsection{Interactions between Heisenberg algebra and $\mathfrak O^k(\g)$}

\begin{theorem} \label{orb:comm}
Let $\g, \alpha \in H^*(X)$. Then we have
\begin{eqnarray*}
 [\Oh(\g), \mathfrak p_{-1}(\alpha)]
 &=&\exp(\hbar \cdot {\rm ad}\,{\mathfrak b}) (\mathfrak p_{-1}(\g \alpha))
 \\
 {[} \mathfrak O^k(\g), \mathfrak p_{-1}(\alpha)] &=& \mathfrak p^{(k)}_{-1}(\g\alpha),\quad k \ge 0.
\end{eqnarray*}
\end{theorem}

\begin{proof}
For simplicity of signs, we assume that the cohomology classes
$\g,\alpha$ are of even degree. It suffices to prove the second
identity.

Recall that $\mathfrak p_{-1}(\alpha)(y) =\frac1{(n-1)!} \sum_{g
\in S_n} \ad g\,(y \otimes \alpha)$, for $y \in
H^*_{\text{orb}}(X^{n-1}/S_{n-1}).$ Regarding $S_{n-1}$ as the
subgroup $S_{n-1}\times 1$ of $S_n$, we introduce an injective
ring homomorphism
$$\iota :H^*(X^{n-1},S_{n-1}) \rightarrow
H^*(X^n,S_n)$$
by sending $\alpha_{\sigma}$ to $\alpha_\sigma \otimes 1_X$, where
$\sigma \in S_{n-1}$. Thus

\begin{eqnarray*}
 && (n-1)! \;{[} \mathfrak O^k(\g), \mathfrak p_{-1}(\alpha)] (y)  \\
 &=& (n-1)! \;[ \mathfrak O^k(\g)\cdot \mathfrak p_{-1}(\alpha)(y)
  - \mathfrak p_{-1}(\alpha)\cdot \mathfrak O^k(\g)(y)]  \\
 &=& O^k(\g,n) \circ \sum_{g \in S_n} \ad g \,(y \otimes \alpha)
   -\sum_{g \in S_n} \ad g \,((O^k(\g,n-1) \circ y) \otimes \alpha)  \\
 &=& \sum_g \ad g \, [(O^k(\g,n) -\iota(O^k(\g,n-1)))\circ (y \otimes \alpha)]
\end{eqnarray*}
where we used the fact that $O^k(\g,n)$ is $S_n$-invariant. By
definition, we have $O^k(\g,n) -\iota(O^k(\g,n-1)) =
(-\xi_{n;n})^{\circ k} \circ \g^{(n)}$. Thus, we obtain
\begin{eqnarray*}
 (n-1)! \;{[} \mathfrak O^k(\g), \mathfrak p_{-1}(\alpha)] (y)
 &=& \sum_g \ad g \,[(-\xi_{n;n})^{\circ k} \circ \g^{(n)} \circ (y \otimes \alpha)]  \\
 &=& \sum_g \ad g \,[  (-\xi_{n;n})^{\circ k}  \circ (y \otimes \g\alpha)].
 \end{eqnarray*}

It remains to prove that
\begin{eqnarray} \label{orb:eqinduction}
\sum_{g\in S_n} \ad g \,[  (-\xi_{n;n})^{\circ k}  \circ (y \otimes \g\alpha)]
=(n-1)! \;\mathfrak p^{(k)}_{-1}(\g\alpha) (y).
\end{eqnarray}

We will prove this by induction. It is clearly true for $k=0$.
Note that $O^1(1_X,n)-\iota( O^1(1_X,n-1)) =-\xi_{n;n}.$ Under the
assumption that the formula (\ref{orb:eqinduction}) is true for
$k$, we have
 \begin{eqnarray*}
  && \sum_g \ad g \,[ (-\xi_{n;n})^{\circ (k+1)} \circ (y \otimes \g\alpha)]\\
 &=&\sum_g \ad g \,[(O^1(1_X,n)-\iota( O^1(1_X,n-1)))
    \circ (-\xi_{n;n})^{\circ k} \circ (y \otimes \g \alpha)]  \\
 &=&O^1(1_X,n) \circ\sum_g \ad g \,[(-\xi_{n;n})^{\circ k} \circ
   (y \otimes \g \alpha)] \\
  && -\sum_g \ad g \,[\iota (O^1(1_X,n-1)) \circ (-\xi_{n;n})^{\circ k}
   \circ (y \otimes \g \alpha)],
 \end{eqnarray*}
since $O^1(\g,n)$ is $S_n$-invariant. By using the induction assumption twice,
we get
  \begin{eqnarray*}
  && \sum_g \ad g \,[ (-\xi_{n;n})^{\circ (k+1)} \circ (y \otimes \g\alpha)]\\
  &=& (n-1)! \; O^1(1_X,n) \circ \mathfrak
  p^{(k)}_{-1}(\g\alpha)(y) \\
   &&-\sum_g \ad g \,[(-\xi_{n;n})^{\circ k} \circ ( (O^1(1_X,n-1) \circ y) \otimes\g \alpha)] \\
  &=& (n-1)!\; [\mathfrak b \cdot \mathfrak p^{(k)}_{-1}(\g\alpha) (y)
   -\mathfrak p^{(k)}_{-1}(\g\alpha) (O^1(1_X,n-1) \circ y) ] \\
   &=& (n-1)! \;\mathfrak p^{(k+1)}_{-1}(\g\alpha) (y) .
 \end{eqnarray*}
So by induction, we have established (\ref{orb:eqinduction}) and
thus the theorem.
\end{proof}

We also have a theorem concerning the operator $\eta(\g)$. It
generalizes Proposition 4.6 in \cite{LS2} (which corresponds to
our special case when $\g =1_X$ and the assumption there that $X$
is a surface is unnecessary).

\begin{theorem}
Let $\g, \alpha \in H^*(X)$ and we further assume that $\g$ can be
written as a sum of classes of even degree. Then we have
\begin{eqnarray*}
 \eta(\g)  \cdot\mathfrak p_{-1}(\alpha)
 &=& \mathfrak p_{-1}(\g\alpha)\cdot \eta(\g) -\mathfrak p_{-1}'(\alpha)\cdot
 \eta(\g), \\
 \varepsilon(\g)  \cdot\mathfrak p_{-1}(\alpha)
 &=& \mathfrak p_{-1}(\g\alpha)\cdot \varepsilon(\g) +\mathfrak p_{-1}'(\alpha)\cdot
 \varepsilon (\g).
\end{eqnarray*}
\end{theorem}

\begin{proof}
The proof of the second formula is similar, so we will prove the
first one only. For simplicity of signs in the proof, we assume
that the cohomology class $\alpha$ is of even degree.

By definition, $\eta_n(\g) =\prod_{i=1}^n (\xi_{i;n} +\g^{(i)})$.
It follows that
$$\eta_n(\g) - \eta_{n-1}(\g) \otimes \g =\xi_{n;n}
\circ \iota(\eta_{n-1}(\g)).$$
Given $y \in H^*_{\text{orb}}(X^{n-1}/S_{n-1})$, we have
 \begin{eqnarray*}
 && (n-1)! \;[ \eta(\g) \mathfrak \cdot \mathfrak p_{-1}(\alpha) (y)
  -\mathfrak p_{-1}(\g \alpha) \cdot \eta(\g) (y)] \\
 &=& \eta_n (\g) \circ \sum_{g \in S_n} \ad g\, (y \otimes \alpha)
   -\sum_{g \in S_n}  \ad g\,[(\eta_{n-1}(\g) \circ y) \otimes \g\alpha] \\
 &=& \sum_g \ad g\,[(\eta_n(\g) - \eta_{n-1}(\g)\otimes \g ) \circ (y \otimes \alpha)]\\
 &=& \sum_g \ad g\,[  \xi_{n;n} \circ\iota(\eta_{n-1}(\g)) \circ (y \otimes \alpha)]  \\
 &=& \sum_g \ad g\,[  \xi_{n;n} \circ ((\eta_{n-1}(\g)  \circ  y) \otimes \alpha)]  \\
 &=& - O^1(1_X,n) \circ \sum_g \ad g\, [(\eta_{n-1}(\g) \circ y) \otimes \alpha]  \\
 && + \sum_g \ad g\, [( O^1(1_X,n-1) \circ \eta_{n-1}(\g) \circ y) \otimes \alpha]  \\
 &=& (n-1)! \;[ -\mathfrak b \cdot \mathfrak p_{-1}(\alpha) \cdot
 \eta(\g) (y) + \mathfrak p_{-1} (\alpha) \cdot \mathfrak b \cdot
 \eta (\g) (y)]  \\
 &=& -(n-1)! \;\mathfrak p_{-1}'(\alpha) \cdot \eta(\g) (y).
 \end{eqnarray*}

This finishes the proof.
\end{proof}
\section{The orbifold cohomology ring of symmetric products  II}
\label{sec:orbformal}

We see from Theorem~\ref{orb:heis}, Theorem~\ref{orb:cubic} and
Theorem~\ref{orb:comm} that the orbifold cohomology rings
$\orbsym$ satisfy the axioms in Subsect.~\ref{subsec:axiom}.
Therefore, we can follow the approaches of [{\bf LQW1-4}] to
establish the results in the following subsections. The terms
involving the canonical class $K$ of $X$ in various formulas in
[{\bf LQW1-4}] will disappear because there is no $K$-term in
Theorem~\ref{orb:cubic}. We remark that the fact $e^2=0$ was used
in \cite{LQW3, LQW4} (where $X$ is a surface). Thus, when dealing
the orbifold cohomology ring $\orbsym$ in this section, we
sometimes need to treat separately and carefully the case when $X$
is a point (i.e. when $e^2 \neq 0$).
\subsection{The ring generators for $\orbsym$}

\begin{theorem} \label{orb:generator}
\begin{enumerate}
\item[{\rm (i)}] Given a closed complex manifold $X$ of even dimension, the orbifold cohomology
ring $\orbsym$ is generated by the cohomology classes
$O^i(\alpha, n)$, where $0 \le i < n$ and $\alpha$ runs over a
fixed linear basis of $H^*(X)$;

\item[{\rm (ii)}] The ring $\orbsym$
is generated by the classes $P_i(\alpha, n)$, where $0\le i<n$ and
$\alpha$ runs over a fixed linear basis of $H^*(X)$.
\end{enumerate}
\end{theorem}

\begin{remark}
Part (i) is the counterpart of Theorem~5.30 in \cite{LQW1}, while
part ~(ii) is the counterpart of Theorem~3.23 in \cite{LQW2}.
\end{remark}
\subsection{The universality of the ring $\orbsym$}

\begin{theorem} \label{orb:univ}
Let $X$ be a closed complex manifold of even dimension. The
orbifold cohomology ring $\orbsym$ is determined uniquely by the
ring $H^*(X)$.
\end{theorem}

We refer to this theorem, Proposition~\ref{orbprop:univ}
and Proposition~\ref{orbprop:univd=0} below
as the universality of the ring $\orbsym$. The theorem follows
from the more quantitative descriptions of the orbifold cup product
of ring generators of $\orbsym$ in Proposition~\ref{orbprop:univ}
and Proposition~\ref{orbprop:univd=0}. It also follows
from combining the results of \cite{LS2} and \cite{FG}.

Let $s \ge 1$, and let $\alpha_1, \ldots, \alpha_s \in H^*(X)$ be
homogeneous cohomology classes. For a partition $\pi =\{ \pi_1,
\ldots, \pi_j \}$ of the set $\{1, \ldots, s \}$, we fix the
orders of the elements listed in each subset $\pi_i$ ($1 \le i \le
j$) once and for all, and define $\ell(\pi) = j$, $\alpha_{\pi_i}
=\displaystyle{\prod_{m \in \pi_i} \alpha_m}$, and the sign ${\rm
sign}(\alpha, \pi)$ by the relation
$$\prod_{i=1}^j \alpha_{\pi_i} = {\rm sign}(\alpha, \pi) \cdot
\prod_{i=1}^s \alpha_i.$$

The choice of the orders for the elements listed in each of the
subsets $\pi_i$, $1 \le i \le \ell(\pi)$ will affect the sign
${\rm sign}(\alpha, \pi)$, but will not affect the long expression
in Proposition~\ref{orbprop:univ} below. We denote by ${\bf
1}_{-k} = \frac{{\mathfrak p}_{-1}(1_X)^k}{k!}$ if $k \ge 0$ and
${\bf 1}_{-k} = 0$ if $k<0.$

\begin{proposition} (Universality) \label{orbprop:univ}
Let $X$ be a closed complex manifold of even dimension $d>0$. Let
$n, s \ge 1$, $k_1, \ldots, k_s \ge 0$, and let $\alpha_1, \ldots,
\alpha_s \in H^*(X)$ be homogeneous. Then, the orbifold product $
O^{k_1}(\alpha_1, n) \circ \cdots \circ O^{k_s}(\alpha_s,n)$ in
$\orbsym$ is a finite linear combination of expressions of the
form

\begin{eqnarray*}
\qquad {\rm sign}(\alpha, \pi) \cdot {\bf 1}_{-\left
(n-\sum\limits_{i=1}^{\ell(\pi)}
  \sum\limits_{j=1}^{m_i-r_i} n_{i, j} \right )}
\prod_{i = 1}^{\ell(\pi)} \left ( \prod_{j = 1}^{m_i-r_i}
\mathfrak p_{-n_{i, j}} \right ) (\tau_{(m_i-r_i)*}(\epsilon_i
\alpha_{\pi_i})) \cdot |0\rangle
\end{eqnarray*}
where $\pi$ runs over all partitions of $\{1, \ldots, s \}$,
$\epsilon_i \in \{1_X, e \}$,
$$r_i = |\epsilon_i|/d \le m_i \le 2+ \sum_{j \in \pi_i} k_j,$$
$0 < n_{i, 1} \le \ldots \le n_{i, m_i-r_i}$,
$\sum\limits_{j=1}^{m_i-r_i} n_{i, j} \le \sum\limits_{j \in
\pi_i} (k_j+1)$ for every $i$, and

\begin{eqnarray*}
\sum_{i = 1}^{\ell(\pi)} \left ( m_i - 2 + \sum_{j = 1}^{m_i-r_i}
n_{i, j} \right ) = \sum_{i=1}^s k_i.
\end{eqnarray*}
Moreover, all the coefficients in this linear combination are
independent of the manifold $X$, the cohomology classes $\alpha_1,
\ldots, \alpha_s$, and the integer $n$.
\end{proposition}

\begin{remark}
This proposition is the counterpart of Proposition~5.1,
\cite{LQW3}.
\end{remark}

For the case $d = 0$ (i.e., $X$ is a point), we adopt the
simplified notations $\mathfrak p_m$ and $O^{k}(n)$ for $\mathfrak
p_m(1_X)$ and $O^{k}(1_X, n)$ respectively. We have the following
analog of Proposition~\ref{orbprop:univ}.

\begin{proposition} \label{orbprop:univd=0}
Let $n, s \ge 1$, $k_1, \ldots, k_s \ge 0$. Then, $O^{k_1}(n)
\circ \cdots \circ O^{k_s}(n)$ is a finite linear combination of
expressions of the form

\begin{eqnarray*}
\qquad {\bf 1}_{-\left (n-\sum\limits_{i=1}^{\ell(\pi)}
  \sum\limits_{j=1}^{m_i-2r_i} n_{i, j} \right )}
\prod_{i = 1}^{\ell(\pi)} \left ( \prod_{j = 1}^{m_i-2r_i}
\mathfrak p_{-n_{i, j}} \right ) \cdot |0\rangle
\end{eqnarray*}
where $\pi$ runs over all partitions of $\{1, \ldots, s \}$, 
$m_i, r_i \in \Z_+$ such that
$2r_i \le m_i \le 2+ \sum_{j \in \pi_i} k_j,$ $0 < n_{i, 1} \le \ldots
\le n_{i, m_i-2r_i}$, $\sum\limits_{j=1}^{m_i-2r_i} n_{i, j} \le
\sum\limits_{j \in \pi_i} (k_j+1)$ for every $i$, and

\begin{eqnarray*}
\sum_{i = 1}^{\ell(\pi)} \left ( m_i - 2 + \sum_{j = 1}^{m_i-2r_i}
n_{i, j} \right ) = \sum_{i=1}^s k_i.
\end{eqnarray*}
Moreover, all the coefficients in this linear combination are
independent of $n$.
\end{proposition}
\subsection{The stability of the ring $\orbsym$}

\begin{theorem} \label{orbth:cupprod}
Let $X$ be a closed complex manifold of even dimension $d$. Let $s
\ge 1$ and $k_i \ge 1$ for $1 \le i \le s$. Fix $n_{i, j} \ge 1$
and $\alpha_{i, j} \in H^*(X)$ for $1 \le j \le k_i$, and fix $n$
with $n \ge \sum\limits_{j=1}^{k_i} n_{i, j}$ for all $1 \le i \le
s$. Then the orbifold cup product
\begin{eqnarray*}
\prod_{i=1}^s \left ( {\bf 1}_{-(n - \sum_{j=1}^{k_i} n_{i, j})}
\prod_{j=1}^{k_i} \mathfrak p_{-n_{i, j}}(\alpha_{i, j}) \cdot
|0\rangle \right )
\end{eqnarray*}
in $\orbsym$ is equal to a finite linear combination of monomials
of the form
\begin{eqnarray*}
 {\bf 1}_{-(n - \sum_{a=1}^N m_{a})}
 \prod_{a=1}^N \mathfrak p_{-m_{a}}(\g_{a}) \cdot |0\rangle
\end{eqnarray*}
where $\sum_{a=1}^N m_a \le \sum_{i=1}^s \sum_{j=1}^{k_i}
n_{i,j}$, and $\g_1, \ldots, \g_N$ depend only on $e,
\alpha_{i,j}$, $1\le i \le s, 1\le j\le k_i$. Moreover, the
coefficients in this linear combination are independent of
$\alpha_{i,j}$ and $n$; they are also independent of $X$ provided
$d> 0$.
\end{theorem}

\begin{remark}
This theorem is the counterpart of Theorem~6.1 in \cite{LQW3}.
\end{remark}
\subsection{The stable ring $\mathfrak R_X$}

Given a finite set $S$ which is a disjoint union of subsets $S_0$
and $S_1$, we denote by ${\mathcal P}(S)$ the set of
partition-valued functions $\rho =(\rho(c))_{c \in S}$ on $S$ such
that for every $c \in S_1$, the partition $\rho(c)$ is required to
be {\it strict} in the sense that $\rho(c) =(1^{m_1(c)} 2^{m_2(c)}
\ldots )$ with $m_r(c) = 0$ or $1$ for all $r \ge 1$.

Now let us take a linear basis $S= S_0 \cup S_1$ of $H^*(X)$ such
that $1_X \in S_0$, $S_0 \subset H^{\rm even}(X)$ and $S_1 \subset
H^{\rm odd}(X)$. If we write $\rho =(\rho (c))_{c \in S}$ and
$\rho(c) =(r^{m_r(c)})_{r \ge 1} =(1^{m_1(c)} 2^{m_2(c)} \ldots
)$, then we introduce the following notations:
\begin{eqnarray*}
\ell(\rho) &=& \sum_{c \in S} \ell(\rho (c))
  = \sum_{c\in S, r\geq 1} m_r(c),  \\
\Vert \rho \Vert &=& \sum_{c \in S} |\rho (c)|
  =\sum_{c\in S, r\geq 1} r \cdot m_r(c),  \\
{\mathcal P}_n(S) &=& \{\rho \in{\mathcal P} (S)\;|\;\; \Vert \rho
\Vert =n\}.
\end{eqnarray*}
Given $\rho=(\rho(c))_{c\in S}=(r^{m_r(c)})_{c\in S,r \ge 1} \in
{\mathcal P}(S)$ and $n \ge 0$, we define
\begin{eqnarray*}
  {\mathfrak p}_{- \rho(c)}(c)
  &=& \prod_{r \ge 1}
  {\mathfrak p}_{-r}(c)^{m_r(c)} = {\mathfrak p}_{-1}(c)^{m_1(c)}
  {\mathfrak p}_{-2}(c)^{m_2(c)} \cdots    \\
 {\mathfrak p}_{\rho}(n)
 &=& {\bf 1}_{-(n-\Vert \rho \Vert)}
  \prod_{c\in S}{\mathfrak p}_{-{\rho}(c)}(c)\cdot\vac \in \orbsym
\end{eqnarray*}
where we fix the order of the elements $c \in S_1$ appearing in
$\displaystyle{\prod_{c \in S}}$ once and for all. It is
understood that ${\mathfrak p}_{\rho}(n)=0$ for $0 \le n < \Vert
\rho \Vert$.

As $\rho$ runs over all partition-valued functions on $S$ with
$\Vert \rho \Vert \le n$, the corresponding ${\mathfrak
p}_{\rho}(n)$ linearly span $\orbsym$, as a corollary to
Theorem~\ref{orb:heis}. According to Theorem~\ref{orbth:cupprod}
(for $s=2$), we can write the orbifold cup product in the ring
$\orbsym$ as
\begin{eqnarray} \label{eq_structure}
{\mathfrak p}_{\rho}(n) \circ {\mathfrak p}_{\sigma}(n) =
\sum_{\nu}d_{\rho\sigma}^\nu {\mathfrak p}_{\nu}(n),
\end{eqnarray}
where $\Vert \nu \Vert \le \Vert \rho \Vert +\Vert \sigma \Vert$,
and the structure coefficients $d_{\rho\sigma}^\nu$ are
independent of $n$. Even though the cohomology classes ${\mathfrak
p}_{\nu}(n)$ with $\Vert \nu \Vert \le n$ in $\orbsym$ are not
linearly independent, we can show that (cf. Lemma~7.1,
\cite{LQW3}) the structure constants $d_{\rho\sigma}^\nu$ in the
formula (\ref{eq_structure}) are uniquely determined from the fact
that they are independent of $n$.

\begin{definition}
The {\it stable ring} associated to a closed complex manifold $X$,
denoted by ${\mathfrak R}_X$, is defined to be the ring with a
linear basis formed by the symbols ${\mathfrak p}_\rho$, $\rho \in
{\mathcal P}(S)$ and with the multiplication defined by
\begin{eqnarray*}
{\mathfrak p}_{\rho} \cdot {\mathfrak p}_{\sigma} = \sum_{\nu}
d_{\rho\sigma}^\nu {\mathfrak p}_{\nu}
\end{eqnarray*}
where the structure constants $d_{\rho\sigma}^\nu$ are from the
relations (\ref{eq_structure}).
\end{definition}

Note that the stable ring does not depend on the choice of a
linear basis $S$ of $H^*(X)$ containing $1_X$ since the operator
${\mathfrak p}_n(\alpha)$ depends on the cohomology class $\alpha
\in H^*(X)$ linearly. Clearly the stable ring ${\mathfrak R}_X$
itself is super-commutative and associative. The ring ${\mathfrak
R}_X$ captures all the information of the orbifold cohomology
rings $\orbsym$ for all $n$, as we easily recover the relations
(\ref{eq_structure}) from the ring ${\mathfrak R}_X$. We summarize
these observations into the following.

\begin{theorem} {\rm (Stability)} \label{orbth:stab}
For a closed complex manifold $X$ of even dimension, the
cohomology rings $\orbsym$, $n\ge 1$ give rise to the stable ring
${\mathfrak R}_X$ which completely encodes the cohomology ring
structure of $\orbsym$ for each $n$. The stable ring ${\mathfrak
R}_X$ depends only on the cohomology ring $H^*(X)$.
\end{theorem}

\begin{remark}
This theorem is the counterpart of Theorem~7.1 \cite{LQW3}. The
stable ring here is the counterpart of the Hilbert ring introduced
in Definition~7.1 of \cite{LQW3}. The stability of the convolution
of symmetric groups (which corresponds to our special case when
$X$ is a point) was due to Kerov and Olshanski (cf. \cite{LT}).
\end{remark}

When $\ell(\rho) = 1$, that is, when the partition $\rho(c)$ is
a one-part partition $(r)$ for some element $c \in S$ and is empty
for all the other elements in $S$, we will simply write
$\mathfrak p_\rho = \mathfrak p_{r,c}$.
Just as in \cite{LQW3}, we can show that the stable ring
$\mathfrak R_X$ is isomorphic to the tensor product $P \otimes E$,
where $P$ is the polynomial algebra generated by $\mathfrak
p_{r,c},\;c\in S_0, r \ge 1$ and $E$ is the exterior algebra
generated by ${\mathfrak p}_{r,c},\;c\in S_1, r \ge 1$.
\subsection{The $\W$ algebras}

In this subsection, we assume that $X$ is a closed complex
manifold of even dimension $d>0$. Results on this section are
counterparts of Sect.~5 of \cite{LQW4}. However, some signs have
been modified due to the sign difference between the two
Heisenberg algebra commutators in the setups of Hilbert schemes
and symmetric products, cf. Theorems~\ref{hilb:heis} and
\ref{orb:heis}. The modification is done by carefully tracing the
procedures in \cite{LQW4}.

Let $\alpha \in H^*(X)$, and $\lambda = (\cdots
(-2)^{m_{-2}}(-1)^{m_{-1}} 1^{m_1}2^{m_2} \cdots)$ be a {\em
generalized partition} of the integer $n = \sum_i i m_i$ whose
part $i\in \Z$ has multiplicity $m_i$. Define $\ell(\lambda) =
\sum_i m_i$, $|\lambda| = \sum_i i m_i = n$, $\lambsq  = \sum_i
i^2 m_i$, $\lambda^! = \prod_i m_i!$, and
\begin{eqnarray*}
\mathfrak p_{\lambda}(\tau_*\alpha) = \left ( \prod_i (\mathfrak
p_i)^{m_i} \right ) (\tau_{\ell(\lambda)*}\alpha).
\end{eqnarray*}
Let $-\lambda$ be the generalized partition whose multiplicity of
$i \in \Z$ is $m_{-i}$. A generalized partition becomes a {\em
partition} in the usual sense if the multiplicity $m_ i = 0$ for
every $i < 0$.

For $p \ge 0$, $n \in \Z$ and $\alpha \in H^*(X)$, define
$\jj^p_n(\alpha) \in \End(\Fock)$ to be
\begin{eqnarray*}
  p! \cdot\left(\sum_{\ell(\lambda) = p+1, |\lambda|=n}
  \frac{1}{\lambda^!} \mathfrak p_{\lambda}(\tau_{*}\alpha)
    + \sum_{\ell(\lambda) = p-1, |\lambda|= n} \frac{\lambsq
     + n^2 - 2}{24\lambda^!} \mathfrak
     p_{\lambda}(\tau_{*}(e\alpha))
\right)
\end{eqnarray*}
where the $\lambda$'s are generalized partitions. Note that
$\jj^0_n(\alpha) = \mathfrak p_{n}(\alpha)$. We define $\Wax$ to
be the linear span of the identity operator ${\rm Id}_{\Fock}$ and
the operators $\jj^p_n(\alpha)$ in $\End(\Fock)$, where $p \ge 0,
n \in \Z$ and $\alpha \in H^*(X)$.

The following theorem describes the operator $\mathfrak
O^k(\alpha)$ in terms of the Heisenberg generators explicitly. It
is a counterpart of Theorem~4.6 in \cite{LQW4}.

\begin{theorem} \label{orb:zeromode}
Let $k \ge 0$, and $\alpha\in H^*(X)$. Then, $\mathfrak O^k(\alpha)
=\frac{(-1)^k}{k+1} \jj^{k+1}_0(\alpha).$
\end{theorem}

In terms of vertex operators, the operator $\jj^p_m(\alpha)$ can
be rewritten as:
\begin{eqnarray} \label{eq:vertex}
 & & \frac{1}{(p+1)} :\mathfrak p^{p+1}:_m (\ta) \nonumber + \frac1{24}
                 p(m^2-3m-2p) :\mathfrak p^{p-1}:_m(\tea)  \\
 & &+ \frac{p(p-1)}{24}  :(\pa^2 \mathfrak p)\, \mathfrak p^{p-2}:_m
 (\tea).
\end{eqnarray}

If we want the coefficients above to be independent of $m$, we can
further rewrite

\begin{eqnarray*}
\jj^p_m(\alpha)
 &=& \frac1{(p+1)} :\mathfrak p^{p+1} :_m (\ta)
 + \frac{p}{24}  (\pa^2: \mathfrak p^{p-1}:)_m(\tea)   \\
 && + \frac{(p+1)p}{12} (\pa :\mathfrak p^{p-1}
 :)_m(\tea)\\
 && + \frac{p(p^2-p-2)}{24}  :\mathfrak p^{p-1} :_m(\tea)   \\
 & &+\frac{p(p-1)}{24}  :(\pa^2 \mathfrak p) \mathfrak p^{p-2}:_m (\tea).
\end{eqnarray*}

The operators $\mathfrak O^p(\alpha), \mathfrak p_{n}(\alpha),$
and $ \jj^p_n(\alpha)$ are related in the following way.
\begin{proposition}
Given $p \ge 0$, $\alpha, \beta \in H^*(X)$, we have
\begin{eqnarray*}
[\mathfrak O^p(\alpha), \mathfrak p_{n}(\beta)] = -n\cdot
\jj^p_n(\alpha\beta).
\end{eqnarray*}
\end{proposition}

We introduce an integer $\Omega_{m,n}^{p,q}$ for $m,n,p,q \in \Z$
as follows:
\begin{eqnarray*}
 \Omega_{m,n}^{p,q}
  &=& mp^3n^2 +3mp^2n^2q -p^2nq +p^2qn^3-3mp^2n^2 +pnq    \\
  && +3m^2pnq -3mpn^2q -m^3q^2p -pqn^3 -mpq +m^3pq    \\
  && +mpq^2 +2mpn^2 -3m^2pnq^2 -2m^2nq +3m^2nq^2 -m^2nq^3 .
\end{eqnarray*}

\begin{theorem} \label{orbth:walg}
Let $X$ be a closed complex manifold of even dimension $d>0$. The
vector space $\Wax$ is closed under the Lie bracket. More
explicitly, for $m,n \in \Z$, and $\alpha, \beta \in H^*(X)$, we
have
\begin{eqnarray*}
 [\jj^p_m(\alpha), \jj^q_n(\beta)]
 &=& (q m -p n) \cdot \jj^{p+q-1}_{m+n} (\alpha\beta)
  + \frac{\Omega_{m,n}^{p,q}}{12} \cdot \jj^{p+q-3}_{m+n}(e\alpha\beta)
\end{eqnarray*}
where $(p,q) \in \Z_+^2$ except for the {\em unordered} pairs
$(0,0), (1,0), (2,0)$ and $(1,1)$. In addition, for these four
exceptional cases, we have
\begin{eqnarray*}
 {[} \jj^0_m(\alpha), \jj^0_n(\beta)]
 &=&  m \delta_{m,-n} \int_X(\alpha\beta) \cdot {\rm Id}_{\Fock}, \\
 {[} \jj^1_m(\alpha), \jj^0_n(\beta)]
 &=& - n \cdot \jj^0_{m+n} (\alpha\beta),\\
 {[} \jj^2_m(\alpha), \jj^0_n(\beta)]
 &=& - 2n \cdot \jj^1_{m+n} (\alpha\beta)
     + \frac{m^3-m}{6}\delta_{m,-n} \int_X(e\alpha\beta)
     \cdot {\rm Id}_{\Fock},\\
 {[} \jj^1_m(\alpha), \jj^1_n(\beta)]
 &=& (m-n) \cdot \jj^1_{m+n} (\alpha\beta)
     + \frac{m^3-m}{12}\delta_{m,-n} \int_X(e\alpha\beta)
     \cdot {\rm Id}_{\Fock}.
\end{eqnarray*}
\end{theorem}

\begin{remark}
This $\W$ algebra should be viewed as a generalization of the
$\W_{\infty}$ algebra (cf. e.g. \cite{FKRW, Kac}). The assumption
$d>0$ above ensures that $e^2 =0$. The case when $d=0$ (i.e. $X$
is a point) has been treated in \cite{LT}.
\end{remark}

\begin{remark} \label{orbrem:zeromode}
Our understanding of $\mathfrak O^k(\alpha)$ in terms of vertex
operators also allows us to recast the study of the ring structure
problems from a different perspective starting from vertex
algebras. Given an integral lattice $A$ (i.e. a free abelian group
with a non-degenerate bilinear form $A \times A \rightarrow \Z$),
the Fock space $\mathfrak F_A$ of a Heisenberg algebra associated
to $A$ affords a natural $\Z_+$-grading: $\mathfrak F_A
=\oplus_{n=0}^\infty \mathfrak F_A^{n}$, and in addition
$\mathfrak F_A$ carries a natural vertex algebra structure
\cite{Bor}. There are numerous operators in $\End (\mathfrak F_A)$
arising from the vertex algebra constructions.

Let us assume that $A$ has an additional structure of a graded
Frobenius algebra compatible with the given bilinear form on $A$.
We may ask if there is any reasonable graded commutative ring
structure on $\mathfrak F_A^{n}$ for each $n$ which comes from the
vertex algebra structure on $\mathfrak F_A$. The answer to this
question is affirmative. We may introduce operators $\mathfrak
O^k(\alpha)$ acting on $\mathfrak F_A^{n}$ for each $n$ to be the
zero-modes of the vertex operators given in
Theorem~\ref{orb:zeromode} and (\ref{eq:vertex}) above. The
operators $\mathfrak O^k(\alpha)$ commute with each other by
Theorem~\ref{orbth:walg}. Next, we define elements $O^k(\alpha,n)$
in $\mathfrak F_A^{n}$ by applying the operator $\mathfrak
O^k(\alpha)$ to $\frac1{n!} \mathfrak p_{-1}(1_A)^n \vac \in
F_A^{n}$, where $ \mathfrak p_{-1}(1_A)$ is a Heisenberg
generator. Then, we define the product in $\mathfrak F_A^{n}$
(which is commutative) by letting $O^a(\alpha,n)\circ O^b(\beta,n)
=\mathfrak O^a(\alpha) \cdot \mathfrak O^b(\beta)\cdot \frac1{n!}
\mathfrak p_{-1}(1_A)^n \vac.$ By Theorem~\ref{orb:generator} we
see that $\mathfrak F_A^{n}$ is generated as a ring by the
elements $O^k(\alpha, n)$'s. In this way, we define a ring
structure on $\mathfrak F_A^{n}$ for each $n$ with a set of ring
generators.
\end{remark}

\section{The deformed orbifold cohomology ring of symmetric products} \label{sec:sign}
\subsection{A deformed orbifold cohomology ring}

Given a complex manifold $M$ with a finite group $G$ action, we
denote by $H^*(M,G;\C)$ and $H^*_{\text{orb}}(M/G;\C)$
respectively the counterparts of $H^*(M,G)$ and
$H^*_{\text{orb}}(M/G)$ with $\C$-coefficients.

\begin{definition}
Let $M$ be a complex manifold with a finite group $G$ action. Let
$t$ be a nonzero complex parameter. We define a product structure,
denoted by $\circ_t$, on $H^*(M,G;\C)$:

$$\alpha_g \circ_t \beta_h = t^{\epsilon(g,h)} \alpha \circ \beta,$$
where $g,h \in G$, $\alpha_g \in H^*(M^g)$, $\beta_h \in
H^*(M^h)$, and $\epsilon (g,h) = (F^g +F^h-F^{gh})/2$. For the
sake of simplicity, we assume here that $(F^g +F^h-F^{gh})/2$ is
an integer for every $g,h\in G$.
\end{definition}

\begin{remark}
In the above definition, for brevity, we have omitted the
dependence of shift numbers on the connected components. If
$\epsilon(g,h)$ is a rational number for some $g,h$, we make sense
of $t^{\epsilon(g,h)}$ by fixing a suitable root of $t$. This
definition is a simple generalization of the signed orbifold
product (i.e. our $t=-1$ case) introduced in \cite{FG}, which in
turn was motivated by Lehn and Sorger \cite{LS2} who introduced
the sign in the symmetric product setup. The new product $\circ_t$
remains to be associative thanks to the identity $\epsilon (g,h)
+\epsilon (gh,k) =\epsilon (g,hk)+\epsilon (h,k)$. It induces a
graded commutative product structure on
$H^*_{\text{orb}}(M/G;\C)$, the $G$-invariant part of
$H^*(M,G;\C)$ (which uses the easy identity $\epsilon(g,h)
=\epsilon(h, h^{ -1}gh)$). Compare with \cite{FG}, Theorem~1.29
and its proof.
\end{remark}

\begin{remark}
In the above definition, if $t^{\epsilon(g,h)}$ are all rational
for a given $t$ and every $g,h\in G$ (e.g. when $t$ is rational),
it makes sense to talk about the ring product $\circ_t$ on the
$H^*(M,G)$ and $H^*_{\text{orb}}(M/G)$ with rational coefficients.
\end{remark}

\begin{proposition} \label{prop:trivialdeform}
The family of ring structures $\circ_t$ on $H^*(M,G;\C)$ (and
resp. on $H^*_{\text{orb}}(M/G;\C)$) with $\C$-coefficient are
isomorphic for all nonzero $t\in \C$.
\end{proposition}

\begin{proof}
For $t \neq 0$, we define the linear map from $H^*(M,G;\C)$ with
product $\circ_t$ to $H^*(M,G;\C)$ with the (original) product
$\circ =\circ_1$:
$$\zeta_t: (H^*(M,G;\C), \circ_t) \rightarrow (H^*(M,G;\C),
\circ)$$ by sending $t^{-F^g/2} \alpha_g$ to $\alpha_g$, for
$\alpha_g \in H^*(M^g)$. The ring isomorphism follows from the
definition of the product $\circ_t$.
\end{proof}

\begin{remark}
The deformed product $\circ_t$, where $t \in \mathbb Q$, on
$H^*(M,G)$ in general can be non-isomorphic for different $t$, and
so this is an interesting deformation. For example, for symmetric
product $X^n/S_n$ associated with a complex surface $X$, the
number $(F^g +F^h-F^{gh})/2$ is always an integer for every
$g,h\in S_n$, but $F^g/2$ often not. So in general $t^{-F^g/2}$
may not be a rational number even when $t$ is, and thus the
isomorphism given in Proposition~\ref{prop:trivialdeform} is not
valid over $\mathbb Q$.

\end{remark}
\subsection{The symmetric product case}
In this subsection, we have two options. Either we assume $t$ is
chosen such that a cubic root $t^{1/3}$ of $t$ is rational, then
all the (orbifold) cohomology groups involved use $\mathbb
Q$-coefficients. On the other hand, if we choose to use
$\C$-coefficients for the cohomology groups, then $t$ can be any
complex number.

Now let $X$ be a closed complex manifold of even dimension $d$.
Let us fix a cubic root $t^{1/3}$ of $t$. We introduce the
modified Heisenberg operators: ${}^t{\mathfrak
p}_n=t^{d/3}\mathfrak p_n$ if $n \le 0$, and ${}^t {\mathfrak
p}_n= t^{-d/6}\mathfrak p_n$ if $n > 0$. Then the Heisenberg
algebra commutation relations in Theorem~\ref{orb:heis} becomes
\begin{eqnarray} \label{orb:deformheis}
 [ {}^t{\mathfrak p}_{m} (\alpha), {}^t{\mathfrak p}_{n}(\beta)]
  = t^{d/6} m \delta_{m,-n} (\alpha, \beta) \cdot \text{Id}_{\Fock}.
\end{eqnarray}
We introduce the modified vertex operator: ${}^t{\mathfrak
p}(\alpha)(z) =\sum_{n\in \Z} {}^t{\mathfrak p}_n(\alpha)
z^{-n-1}.$

We modify the definition of the operators ${\mathfrak O}(\alpha)$,
${\mathfrak O}^k(\alpha)$, and ${\mathfrak O}_{\hbar}(\alpha)$ by
using the product $\circ_t$ instead of $\circ$, and denote by the
resulting operators by ${}^t{\mathfrak O}(\alpha)$,
${}^t{\mathfrak O}^k(\alpha)$, and ${}^t{\mathfrak O}_\hbar
(\alpha)$. We denote by ${}^t{\mathfrak b} ={}^t{\mathfrak
O}^1(1_X).$ Given $\mathfrak f \in \End(\Fock)$, we denote by
$({\rm ad}\,{}^t{\mathfrak b} )\mathfrak f  =[ {}^t{\mathfrak b},
\mathfrak f].$

The same argument as earlier leads to the following theorem.

\begin{theorem} \label{orb:deformcomm}
Let $\g, \alpha \in H^*(X)$. Then we have
\begin{eqnarray*}
 [{}^t{\mathfrak O}_\hbar(\g), {}^t{\mathfrak p}_{-1}(\alpha)]
 &=&\exp(\hbar \cdot {\rm ad}\,{}^t{\mathfrak b}) ({}^t{\mathfrak p}_{-1}(\g \alpha))
 \\
 {[} {}^t{\mathfrak O}^k(\g), {}^t{\mathfrak p}_{-1}(\alpha)]
 &=& ({\rm ad}\,{}^t{\mathfrak b})^k \cdot {}^t{\mathfrak p}_{-1}(\g\alpha),\quad k \ge 0.
\end{eqnarray*}
\end{theorem}

\begin{theorem} \label{orb:deformcubic}
We have ${}^t{\mathfrak b} =-\frac16 :{}^t{\mathfrak p}^3:_0
(\tau_{*}1_X).$
\end{theorem}

\begin{remark}   \label{orb:deformproofcubic}
An outline of a proof of Theorem~\ref{orb:deformcubic} goes as
follows. For notational simplicity, we suppress the dependence on
cohomology classes of $\mathfrak p_n$, $:{}^t{\mathfrak p}^3:_0$
below. Given $\sigma \in S_n$, it is well known that the degree
shift number $F^\sigma =\frac{d}2 \cdot d(\sigma)$. Let us look at
the product of a transposition $(a,b)$ with $\alpha \in
H^*((X^n)^\sigma)$. If $a, b$ do not lie in the same cycle of
$\sigma$, then $(a,b)\sigma$ is obtained from $\sigma$ by
combining the two cycles in $\sigma$ containing $a,b$
respectively. It follows that $\epsilon ((a,b), \sigma) =0$. If
$a,b$ lies in a same cycle of $\sigma$, then $(a,b)\sigma$ is
obtained from splitting the cycle of $\sigma$ containing $a$ and
$b$ into two cycles, and thus $\epsilon ((a,b), \sigma)
=\frac{d}2$. See the proof of Theorem~2, \cite{FW}, for some
illustration by concrete examples. This explains the factor
$t^{d/2}$ below (compare with Remark~\ref{orb:proofcubic}):

\begin{eqnarray*} \label{orbeq:modifycubic}
{}^t{\mathfrak b} =-\hf \sum_{n,m >0} (\mathfrak p_{-n-m}
\mathfrak p_{n} \mathfrak p_{m} + t^{d/2} \mathfrak p_{-n}
\mathfrak p_{-m} \mathfrak p_{n+m}) .
\end{eqnarray*}
On the other hand, one observes that

$$-\frac16
:{}^t{\mathfrak p}^3:_0
 = -\hf \sum_{n,m>0} ({}^t\mathfrak p_{-n-m}
{}^t\mathfrak p_{n} {}^t\mathfrak p_{m} + {}^t\mathfrak p_{-n}
{}^t\mathfrak p_{-m} {}^t\mathfrak p_{n+m}), $$
which coincide with ${}^t{\mathfrak b}$ by using the definition of
${}^t\mathfrak p$ and the above formula for ${}^t{\mathfrak b}$.
\end{remark}

\begin{remark}
Formula (\ref{orb:deformheis}), Theorem~\ref{orb:deformcomm} and
Theorem~\ref{orb:deformcubic} indicate that the counterparts of
Theorem~\ref{orb:heis}, Theorem~\ref{orb:cubic} and
Theorem~\ref{orb:comm} hold if we use the product $\circ_t$ on
$\orbsym$ instead of $\circ$. Therefore the results established in
section~\ref{sec:orbformal} also carry over for the product
$\circ_t$ with appropriate modifications.
\end{remark}

\begin{remark}
Starting from a graded Frobenius algebra $A$, we can obtain a
family of Frobenius algebra structures on $\mathfrak F_A^n$
depending on $t$ by using the modified Heisenberg algebra etc.
When setting $t=-1$, the algebra $\mathfrak F_A^n$ should be
isomorphic to $A^{[n]}$ given in \cite{LS2}.
\end{remark}
\subsection{A cohomology ring isomorphism}
Let $X$ be a projective surface. We have seen that both $\Fock
=\oplus_n \orbsym$ and $\Fockh =\oplus_n H^*(\Xn)$ are Fock spaces
of the same size. By sending $\mathfrak p_{-n_1}(\alpha_1)\cdots
\mathfrak p_{-n_k}(\alpha_k)\vac$ to $\mathfrak
a_{-n_1}(\alpha_1)\cdots \mathfrak a_{-n_k}(\alpha_k)\vac$, where
$n_1, \ldots, n_k >0$, $\alpha_1, \ldots, \alpha_k \in H^*(X)$, we
have defined a linear isomorphism $\Theta: \Fock \rightarrow
\Fockh$, which induces a linear isomorphism $\Theta_n: \orbsym
\rightarrow H^*(\Xn)$ for each $n$.

\begin{theorem}
Let $X$ be a projective surface with numerically trivial canonical
class. The linear map $\Theta_n: \orbsym \rightarrow H^*(\Xn)$ is
a ring isomorphism, if we use the product $\circ_{-1}$ on
$\orbsym$.
\end{theorem}

\begin{proof}
Noting that $d =2$ and $t =-1$, we can take $t^{d/6}=-1$. Thus,
${}^t{\mathfrak p}_n= \mathfrak p_n$ if $n \le 0$, and ${}^t
{\mathfrak p}_n= -\mathfrak p_n$ if $n > 0$ (we keep using $t$
instead of $-1$ here and below for notational convenience.) The
Heisenberg algebra commutators for the ${}^t {\mathfrak p}_n$ in
(\ref{orb:deformheis}) and for the $\mathfrak a_n$ in
Theorem~\ref{hilb:heis} exactly match. Comparing
(\ref{eq_general}) and (\ref{eq:expon}), we see that $\Theta$
sends $G^1(1_X,n) =c_1({\mathcal O}^{[n]})$ (where $\mathcal O$
denotes the trivial line bundle over $X$) to
$O^1(1_X,n)=-\sum_{i=1}^n \xi_i$. Note that the operator
$\mathfrak d$ and ${}^t\mathfrak b$ are defined in terms of
$G^1(1_X,n)$ and $O^1(1_X,n)$ respectively. By
Theorem~\ref{hilb:cubic} and Theorem~\ref{orb:deformcubic}, the
operator $\mathfrak d$ matches exactly with ${}^t\mathfrak b$.
Then it follows from comparing Theorems~\ref{hilb:comm} and
Proposition~\ref{orb:deformcomm} that the operator $\mathfrak
G^k(\g)$ coincide with $\mathfrak O^k(\g)$. If we recall the
definitions of $\mathfrak G^k(\g)$ and $\mathfrak O^k(\g)$, the
theorem follows now from Theorem~\ref{orb:generator}~(i) and its
Hilbert scheme counterpart Theorem~1.2 in \cite{LQW1}.
\end{proof}

\begin{remark}
This ring isomorphism has been earlier established in a different
way by combining the results in \cite{LS2, FG} (also cf.
\cite{Uri}).
\end{remark}

Modifying $\Theta$, we introduce a linear isomorphism
$\widetilde{\Theta}: \Fock \rightarrow \Fockh$ by sending
$\sqrt{-1}^{\sum_{a=1}^k n_a-k} \mathfrak p_{-n_1}(\alpha_1)\cdots
\mathfrak p_{-n_k}(\alpha_k)\vac$ to $\mathfrak
a_{-n_1}(\alpha_1)\cdots \mathfrak a_{-n_k}(\alpha_k)\vac$. This
induces a linear isomorphism $\widetilde{\Theta}_n:
H^*_{\text{orb}}(X^n/S_n;\C) \rightarrow H^*(\Xn;\C)$ for each
$n$. Note that both $\Theta_1$ and $\widetilde{\Theta}_1$ are
simply the identity map on the cohomology group of the surface
$X$.

\begin{theorem}
Let $X$ be a projective surface with numerically trivial canonical
class. The linear map $\widetilde{\Theta}_n:
H^*_{\text{orb}}(X^n/S_n;\C) \rightarrow H^*(\Xn;\C)$ is a ring
isomorphism from the cohomology ring of Hilbert scheme with
$\C$-coefficient to the standard orbifold cohomology ring of the
symmetric product with $\C$-coefficient.
\end{theorem}

\begin{proof}
Note that $\mathfrak p_{-n}(\alpha) \vac$ corresponds to an
$n$-cycle whose shift number is $(n-1)$ and a permutation
associated to $\mathfrak p_{-n_1}(\alpha_1)\cdots \mathfrak
p_{-n_k}(\alpha_k)\vac$ has shift number $\sum_{a=1}^k n_a-k$.
Thus, the map $\widetilde{\Theta}_n$ is the composition of the
ring isomorphism $\Theta_n$ with the ring isomorphism $\zeta_t$
for $t=-1$ defined in the proof of
Proposition~\ref{prop:trivialdeform}.
\end{proof}

The above theorem supports the original conjecture of Ruan
\cite{Ru1} if the cohomology coefficient is $\C$ rather than
$\mathbb Q$. Of course, the surface example at the end of
Section~2 of \cite{FG} is no longer a counterexample over $\C$,
since all symmetric bilinear form over $\C$ can be diagonalizable.
Our results refresh the hope that Ruan's Conjecture may be valid
for any hyperkahler resolution, once we insist on the cohomology
coefficient being $\C$.
\section{Open questions} \label{sec:open}

In this section, we list some open problems for further research.

{\bf Question $1$.} Understand the cohomology ring structure of
the Hilbert scheme $\Xn$ when $X$ is a quasi-projective surface.
While certain degeneracies occur in connection between the Chern
character operators and vertex operators (cf. \cite{Lehn},
Sect.~4.4, and \cite{LS1} for the affine plane case), we expect
that most of the geometric statements, such as those on ring
generators, universality and stability, should remain valid in the
quasi-projective case.

{\bf Question $1'$.} Understand the orbifold cohomology ring
structure of the symmetric products $X^n/S_n$ for a non-closed
complex manifold $X$.

{\bf Question $2$.} Use the axiomatization in
Sect.~\ref{subsec:axiom} to check Ruan's conjecture on the
isomorphism between the (signed) orbifold cohomology ring of the
symmetric product $X^n/S_n$ and the quantum corrected cohomology
ring of the Hilbert scheme $\Xn$, when $X$ is an arbitrary
(quasi-)projective surface.

{\bf Question $3$.} Is there a family of ring structures on the
rational cohomology group of the Hilbert scheme $\Xn$ depending on
a rational parameter $t$, such that when $t=-1$ it is the standard
one and that it becomes isomorphic to the deformed orbifold
cohomology ring $(\orbsym, \circ_t)$ when $X$ has a numerically
trivial canonical class? We may ask similar questions for crepant
resolutions of orbifolds.

{\bf Question $4$.} Why is the theory of vertex algebras so
effective in the study of the geometry of Hilbert schemes and
symmetric products? On the other hand, when the canonical class
$K$ of the surface $X$ is not numerically trivial, $K$ becomes an
obstruction in connection between Hilbert schemes and vertex
algebras. How is this related to the quantum corrections on
Hilbert schemes as proposed by Ruan?

{\bf Question $5$.} The appearance of $\mathcal W$ algebras
indicates connections to completely integrable systems. How to see
this in the framework of Hilbert schemes and symmetric products?

{\bf Question $6$.} How to understand the orbifold cohomology ring
of the symmetric products $X^n/S_n$ for $X$ of odd complex
dimension, or even for a more general manifold $X$?

We end this paper with a table comparing the pictures of Hilbert
schemes and symmetric products (see above). The reader may compare
with another table in \cite{Wa2} which relate the pictures of
Hilbert schemes and wreath products.

\begin{table}
 \begin{center}
\caption{A DICTIONARY} \label{Tab_dict}
\begin{tabular}{|| l | l || }
\hline
 Hilbert Scheme $\Xn$         &  Symmetric Product $X^n/S_n$  \\ \hline\hline
 $\Fockh=\oplus_nH^*(\Xn)$    &  $\Fock =\oplus_n \orbsym$   \\ \hline
 cup product                  &  (signed) orbifold cup product   \\ \hline
 Heisenberg generator $\mathfrak a_n(\alpha)$
                 &   Heisenberg generator $\mathfrak p_n(\alpha)$  \\ \hline
 total Chern class $c(L^{[n]})$
                              &  $\varepsilon_n(c(L))$  \\ \hline
 $c(L^{[n]\vee})$
                              & $\eta_n(c(L^\vee))$  \\ \hline
 Chern roots of $L^{[n]}$
                              & $ c_1(L)^{(i)}-\xi_i $ \\ \hline
 class $G^k(\alpha,n)$        &   class $O^k(\alpha,n)$    \\ \hline
 Lehn's operator $\mathfrak d$
                              &   generalized Goulden's operator $\mathfrak b$    \\ \hline
\end{tabular}
 \end{center}
\end{table}
\vspace{1cm}


\begin{thebibliography}{ABCD}

\bibitem[BBM]{BBM} P. Baum, J. Brylinski and R. MacPherson,
{\em Cohomologie \'equivariante d\'elocalis\'ee}, C.R. Acad. Sci.
Paris {\bf 300} (1985), 605--608.

\bibitem[BC]{BC} P. Baum and A. Connes,
{\em Chern character for discrete groups}, In: Y.~Matsumoto et al
(eds.), A Fete of Topology, Academic Press, 1988.

\bibitem[Bor]{Bor} R.~Borcherds,
{\em Vertex algebras, Kac-Moody algebras, and the Monster}, Proc.
Natl. Acad. Sci. USA {\bf 83} (1986), 3068--3071.

\bibitem[CR]{CR}
W. Chen and Y.~Ruan, {\em A new cohomology theory for orbifold},
math.AG/0004129.

\bibitem[DHVW]{DHVW}
L. Dixon, J.A. Harvey, C. Vafa, and E. Witten, {\em Strings on
orbifolds}, Nuclear Phys. {\bf B 261} (1985), 678--686.

\bibitem[FG]{FG} B. Fantechi, L. G\" ottsche,
{\em Orbifold cohomology for global quotients},
math.AG/0104207.

\bibitem[FKRW]{FKRW} E.~Frenkel, V.~Kac, A. Radul, and W. Wang,
{\em $\W_{1+\infty}$ and ${\mathcal W}(gl_N)$ with central
charge~$N$}, Commun. Math. Phys. {\bf 170} (1995), 337--357.

\bibitem[FW]{FW} I. Frenkel and W. Wang,
{\em Virasoro algebra and wreath product convolution}, J.~Alg.
{\bf 242} (2001), 656--671.

\bibitem[Gou]{Gou} I.~Goulden,
{\em A differential operator for symmetric functions and the
combinatorics of multiplying transpositions}, Trans. Amer. Math.
Soc. {\bf 344} (1994), 421--440.

\bibitem[Gro]{Gro} I.~Grojnowski,
{\em Instantons and affine algebras I: the Hilbert scheme and
vertex operators}, Math. Res. Lett. {\bf 3} (1996), 275--291.

\bibitem[Juc]{Juc} A. Jucys,
{\em Symmetric polynomials and the center of the symmetric group
rings}, Rep. Math. Phys. {\bf 5} (1974), 107--112.

\bibitem[Kac]{Kac} V. Kac,
{\em Vertex Algebras for Beginners}, Second Edition, University
Lecture Series {\bf 10}, AMS, Providence, Rhode Island, 1998.

\bibitem[Kuhn]{Kuhn} N.~Kuhn,
{\em Character rings in algebraic topology}, In: {Advances in
Homotopy}, London Math. Soc. Lect. Notes Series {\bf 139} (1989),
111--126.

\bibitem[LT]{LT} A. Lascoux and J.-Y. Thibon,
{\em Vertex operators and the class algebras of symmetric groups},
Preprint, math.CO/0102041.

\bibitem[Lehn]{Lehn} M. Lehn,
{\em Chern classes of tautological sheaves on Hilbert schemes of
points on surfaces}, Invent. Math. {\bf 136} (1999), 157--207.

\bibitem[LS1]{LS1} M. Lehn and C. Sorger,
{\em Symmetric groups and the cup product on the cohomology of
Hilbert schemes}, Duke Math. J. (to appear), math.AG/0009131.

\bibitem[LS2]{LS2} ------,
{\em The cup product of the Hilbert scheme for $K3$ surfaces},
math.AG/0012166.

\bibitem[LQ]{LQ} W.-P. Li and Z. Qin,
{\em On $1$-point Gromov-Witten invariants of the Hilbert schemes
of points on surfaces}, Preprint.

\bibitem[LQW1]{LQW1} W.-P. Li, Z. Qin, and W. Wang, {\em Vertex algebras and the
cohomology ring structure of Hilbert schemes of points on
surfaces}, Math. Ann. (to appear), math.AG/0009132.

\bibitem[LQW2]{LQW2} ------, {\em Generators
for the cohomology ring of Hilbert schemes of points on surfaces},
Intern. Math. Res. Notices No. {\bf 20} (2001) 1057--1074,
math.AG/0009167.

\bibitem[LQW3]{LQW3} ------, {\em Universality and
stability of cohomology rings of Hilbert schemes of points on
surfaces}, Preprint, math.AG/0107139.

\bibitem[LQW4]{LQW4} ------,
{\em Hilbert schemes and $\mathcal W$ algebras}, 
Intern. Math. Res. Notices (to appear),
math.AG/0111047.

\bibitem[Mac]{Mac} I.~G. Macdonald,
{\em Symmetric functions and Hall polynomials}, 2nd Ed., Clarendon
Press, Oxford, 1995.

\bibitem[Mur]{Mur} G. Murphy,
{\em A new construction of Young's seminormal representation of
the symmetric group}, J. Alg. {\bf 69} (1981), 287--291.

\bibitem[Na1]{Na1} H. Nakajima,
{\em Heisenberg algebra and Hilbert schemes of points on
projective surfaces}, Ann. Math. {\bf 145} (1997), 379--388.

\bibitem[Na2]{Na2} ------,
{\em Lectures on Hilbert schemes of points on surfaces}, Univ.
Lect. Ser. {\bf 18}, Amer. Math. Soc. (1999).

\bibitem[Ru1]{Ru1} Y.~Ruan,
{\em Stringy geometry and topology of orbifolds}, math.AG/0011149.

\bibitem[Ru2]{Ru2} ------,
{\em Cohomology ring of crepant resolutions of orbifolds},
math.AG/0108195.

\bibitem[Seg]{Seg} G.~Segal,
{\em Equivariant K-theory and symmetric products}, Preprint, 1996.

\bibitem[Uri]{Uri} B.~Uribe,
{\em Orbifold Cohomology of the Symmetric Product},
math.AT/0109125.

\bibitem[VW]{VW} C. Vafa and E. Witten,
{\em A strong coupling test of $S$-duality}, Nucl. Phys. {\bf B
431} (1994), 3--77.

\bibitem[Wa1]{Wa1} W.~Wang,
{\em Equivariant K-theory, wreath products, and Heisenberg
algebra}, Duke Math. J. {\bf 103} (2000), 1--23.

\bibitem[Wa2]{Wa2} ------,
{\em Algebraic structures behind Hilbert schemes and wreath
products}, Contemp. Math. (to appear), math.QA/0011103.

\bibitem[Zas]{Zas} E. Zaslow,
{\em Topological orbifold models and quantum cohomology rings},
Commun. Math. Phys. {\bf 156} (1993), 301--331.

\end{thebibliography}
\end{document}